\documentclass[12pt,a4paper]{article}
\usepackage[cp1251]{inputenc} 
\usepackage[russian]{babel}
\usepackage{amsfonts}

\newcommand{\di}{\displaystyle}
\newcommand{\B}{$\hfill\Box$}
\newcommand{\al}{\alpha}
\newcommand{\be}{\beta}
\newcommand{\ga}{\gamma}
\newcommand{\de}{\delta}
\newcommand{\la}{\lambda}
\newcommand{\om}{\omega}
\newcommand{\ee}{\varepsilon}
\newcommand{\vv}{\varphi}
\newcommand{\iy}{\infty}

\begin{document}

\begin{center}
{\large\bf
Recovering Differential Operators on Spatial Networks}\\[0.1cm]
{\bf V.A.\,Yurko} \\[0.1cm]
\end{center}

\thispagestyle{empty}

{\bf Abstract.} We give a short review of results on inverse spectral problems for
ordinary differen\-tial operators on a spatial networks (geometrical graphs).
We pay the main attention to the most important nonlinear inverse problems of
recovering coefficients of differential equations from spectral characteristics
provided that the structure of the graph is known a priori. In the first half of
the review we provide results related to inverse Sturm-Liouville problems on
arbitrary compact graphs. Further, results on inverse problems for arbitrary order
differential operators on compact graphs are presented. At the end we provide
the main results on inverse problems on noncompact graphs.

Bibliography: 55 titles.



\begin{center}
{\large Contents}
\end{center}

\noindent 1. Introduction\\
2. Sturm-Liouville operators on trees.

   2.1. Weyl functions. Statements of the inverse problems.

   2.2. Properties of spectral characteristics.

   2.3. Solution of the inverse problems.\\
3. Sturm-Liouville operators on compact graphs with cycles.

   3.1. Inverse problems on the bush-graph.

   3.2.  Inverse problems on an arbitrary compact graph.\\
4. Arbitrary order differential operators on compact graphs.

   4.1. Inverse problems for $n$-th order operators.

   4.2. Inverse problems for variable order operators on graphs without cycles.

   4.3. Inverse problems for variable order operators on graphs with cycles.\\
5. Differential operators on noncompact graphs.

   5.1. Inverse Sturm-Liouville problems on a noncompact star-graph.

   5.2. Inverse Sturm-Liouville problems on a noncompact sun-graph.

   5.3. Inverse Sturm-Liouville problems on a noncompact A-graph.

   5.4. Variable order operators on a noncompact with a cycle.\\
6. Bibliography.

\begin{center}
{\bf 1. Introduction}
\end{center}

This paper presents a short review of results on inverse spectral problems for ordinary
differential operators on spatial networks (geometrical graphs). Differential operators
on graphs often appear in mathematics, mechanics, physics, geophysics, physical
chemistry, electronics, nanotechnology and other branches of natural sciences and
engineering (see [1-8] and the bibliog\-raphy therein). Most of the works are devoted
to the so-called direct problems of studying properties of the spectrum and root
functions. In this direction we mark the monograph [8], where extensive references
can be found related to direct problems of spectral analysis on graphs. Inverse
problems consist in recovering coefficients of operators from their spectral
characteristics. Such problems have many applications and therefore attract
attention of scientists. On the other hand, inverse spectral problems
are very difficult for investigation, because of their nonlinearity. The main
results and methods on inverse spectral problems for differential operators
{\it on an interval} (finite or infinite) were obtained in the second half of
the XX-th century; these results are presented fairly completely in the monographs
[9-12], where an extensive bibliography can be found.

Active Research on the inverse spectral problems for
differential operators {\it on graphs}  started in the XXI-th centure.
The first important and rather wide class of problems, for
which the theory of solutions of inverse problems has been constructed, were the
so-called {\it trees}, i.e. graphs without cycles [13-15]; see Section 2 for more
details on these results. Two different methods were applied in the investigation
of inverse problems on trees in [13-15]: the transformation operator method and
the method of spectral mappings. The first method [13] connects with the process of
wave propagation, and the second one [14] uses ideas of contour integration and the
apparatus of the theory of analytic functions. In particular, the method of spectral
mappings allowed one not only to prove a uniqueness theorem with most limited data,
but also to obtain an algorithm for the solution of this class of inverse problems.
Moreover, this method turned out to be effective also for essentially more difficult
problems on graphs with cycles. Inverse spectral problems for Sturm-Liouville
operators on arbitrary compact graphs with cycles have been solved in [16-20] and
other works (see also Section 3 of the present paper). Some other aspects of the
inverse problem theory for Sturm-Liouville operators on compact graphs were discussed
in [21-29] and other works. In particular, a generalization of results from [14]
was obtained in [21] for Sturm-Liouville operators on a tree with generalized
matching conditions and with a complex-valued potential; an algorithm for the
solution of the inverse problem from the given Weyl vector was provided, and the
uniqueness of the solution was proved. In [22], for a particular case of the problem
from [21], the uniqueness theorem was proved, which coincides with the uniqueness
theorem from [14], and also the calculation of the lengths of the graph's edges
from the given matrix of spectral mappings was discussed. In [25] a compact graph
with one cycle and one boundary edge is investigated. In [24] the case of singular
potentials from $W_2^{-1}$ was considered, and in [27-29] inverse problems with
incomplete information were studied together with inverse nodal problems on graphs.
In [30-31] inverse problems on graphs for pencils of differential operators
with nonlinear dependence on the spectral parameter were investigated.

Inverse problems {\it on noncompact graphs} are not as thoroughly investigated.
Several works on the topic exist so far, but they do not form a general picture.
For some particular cases of noncompact
graphs, the inverse Sturm-Liouville problem was investigated in [32-43] and other
works. We mark results from [37-40], where the solutions of the inverse problems
were obtained for some important classes of noncompact graphs both with cycles and
without them. These results are represented in Section 5 of the present paper.
However, nowadays the theory is absent for noncompact graphs of an arbitrary structure.

All above-mentioned results are related only to the Sturm-Liouville operator.
Inverse problems for arbitrary order differential operators on graphs, because of
their complexity, have been investigated very little. For $n$-th order operators
for $n>2,$ results are known only for trees. An important and interesting class is
variable order operators, when differential equations have different orders on
different edges. Such equations appear in various problems of mathematics and in
applications (see, for example, [47] and the bibliography therein), for example,
for describing oscillations in such structures as masts with cable supports,
cable-stayed bridges and others. Some results on inverse problems for variable
order differential operators were obtained in [48-52]. These results are
presented in subsections 4.2, 4.3, 5.4 of the present review.

Some words about the structure of the paper.
In Section 2 we provide results related to inverse Sturm-Liouville on compact trees.
Sturm-Liouville operators on compact graphs with cycles are studied in Section 3.
Section 4 is devoted to inverse problems for arbitrary order differential operators
on compact graphs. In Section 5 we present the main results on noncompact graphs.
We note that the review is not meant to be comprehensive. Here we describe only main
results on the theory of inverse spectral problems for ordinary differential operators
on spatial networks. We pay the main attention to the most important nonlinear
inverse problems of recovering coefficients of differential equations from spectral
characteristics provided that the structure of the graph is known a priori. We
mention only the most important monographs and papers, and refer the reader to them
for more detailed acquaintance with the subject.

\begin{center}
{\bf 2. Sturm-Liouville operators on trees}
\end{center}

In this section we provide results on inverse spectral problems for self-adjoint and
non-self-adjoint Sturm-Liouville operators on compact graphs without cycles (trees)
with generalized matching conditions in interior vertices and boundary conditions in
boundary vertices. We introduce and study the so-called Weyl vector which is a
generalization of the notion of the Weyl function for the classical Sturm-Liouville
operator on an interval. We consider three inverse problems of recovering the
potential and the coefficients of the boundary conditions from the Weyl vector,
froma system of spectra and from discrete spectral data, respectively. For these
inverse problems uniqueness theorems are proved, and algorithms for constructing
solutions are provided. For studying these inverse problems we use ideas of the
method of specral mappings, and also an important procedure of "pseudo-cutting"
the tree, which allows one to solve the inverse problems step by step going on
from one edge to another.

{\bf 2.1. Weyl functions. Statements of the inverse problems.}
Consider a compact connected tree $T$ in ${\bf R^m}$ with the root $v_0$, the set
of vertices $V=\{v_0,\ldots,v_r\}$ and the set of edges ${\cal E}=\{e_1,\ldots, e_r\}$.
For simplicity we suppose that the length of each edge is equal to $1$
(the results remain true for arbitrary lengths).
A vertex is called a boundary vertex if it belongs to only one edge.
Such an edge is called a boundary edge. All other vertices and edges are
called internal. Without loss of generality we assume that $v_0$ is a
boundary vertex.

For two points $a,b\in T$ we will write $a\le b$ if $a$ lies on a unique
simple path connecting the root $v_0$ with $b$; let $|b|$ stand for
the length of this path. We will write $a<b$ if $a\le b$ and $a\ne b.$
The relation $<$ defines a partial ordering on $T.$
If $a<b$ we denote $[a,b]:=\{z\in T:\; a\le z\le b\}.$ In particular,
if $e=[v,w]$ is an edge, we call $v$ its initial point, $w$ its end
point and say that $e$ emanates from $v$ and terminates at $w.$ \\ \\ \\

\begin{picture}(100,130)
\unitlength=1.4mm

\put(55,22){\line(6,-1){20}}

\put(55,22){\circle*{0.5}}

\put(61,4.6){\circle*{0.5}}

\put(75,18.65){\line(1,0){20}}

\put(75,18.65){\circle*{0.5}}

\put(95,18.65){\circle*{0.5}}

\put(55,22){\line(-3,1){20}}

\put(55,22){\circle*{0.5}}

\put(16.7,21.4){\circle*{0.5}}

\put(55,22){\line(-2,-1){20}}

\put(75,18.65){\line(-1,-1){14}}

\put(75,18.65){\circle*{0.5}}

\put(35,28.7){\line(-5,3){18}}

\put(35,28.7){\circle*{0.5}}

\put(17,39.5){\circle*{0.5}}

\put(35,28.7){\line(-5,-2){18}}

\put(35,12){\line(-6,1){18}}

\put(35,12){\circle*{0.5}}

\put(17,15){\circle*{0.5}}

\put(17.7,1.6){\circle*{0.5}}

\put(35,12){\line(-5,-3){17.5}}

\put (34,30){$v_8$}

\put (43.5,26.5){$e_8$}

\put (73,20){$v_6$}

\put (54,23){$v_7$}

\put (64,22){$e_7$}

\put (44,15){$e_9$}

\put (35,10){$v_9$}

\put (25,35.5){$e_5$}

\put (13,39){$v_5$}

\put (23,26.5){$e_4$}

\put (13,21){$v_4$}

\put (24,14.8){$e_3$}

\put (12.5,14.5){$v_3$}

\put (25.5,5){$e_2$}

\put (14,0){$v_2$}

\put (68,10){$e_1$}

\put (58,2.5){$v_1$}

\put (83,20){$e_6$}

\put (95.5,18){$v_0$}

\end{picture}

\centerline{\small fig. 2.1}

\bigskip
For each internal vertex $v$ we denote by $R(v):=\{e\in {\cal E}:\;
e=[v,w],\, w\in V\}$ the set of edges emanating from $v.$
For any $v\in V$ the number $|v|$ is a non-negative integer, which is
called the order of $v.$ For $e\in {\cal E}$ its order is defined as
the order of its end point. The number $\sigma:=\max_{j=\overline{1,r}}
|v_j|$ is called the height of the tree $T.$ Let $V^{(\mu)}:=
\{v\in V:\; |v|=\mu\}$, $\mu=\overline{0,\sigma}$ be the set of
vertices of order $\mu,$ and let ${\cal E}^{(\mu)}:=\{e\in {\cal E}:\;
e=[v,w],\, v\in V^{(\mu-1)}, w\in V^{(\mu)}\}$, $\mu=\overline{1,\sigma}$
be the set of edges of order $\mu.$

Each edge $e\in {\cal E}$ is
parameterized by the parameter $x\in [0,1].$ It is convenient for us
to choose the following orientation on each edge $e=[v,w]\in {\cal E}$:
if $z=z(x)\in e,$ then $z(0)=w,\, z(1)=v,$ i.e. $x=0$ corresponds to
the end point $w,$ and $x=1$ corresponds to the initial point $v.$
For definiteness we enumerate the vertices $v_j$ as follows:
$\Gamma:=\{v_0,v_1,\ldots, v_p\}$ are the boundary vertices,
$v_{p+1}\in V^{(1)}$, and $v_j,\, j>p+1$ are enumerated in order of
increasing $|v_j|$. We enumerate the edges similarly, namely:
$e_j=[v_{j_k},v_j]$, $j=\overline{1,r},\, j_k<j.$ In particular,
$E:=\{e_1,\ldots, e_{p+1}\}$ is the set of boundary edges, $e_{p+1}=
[v_0,v_{p+1}].$ Clearly, $e_j\in{\cal E}^{(\mu)}$ iff $v_j\in V^{(\mu)}$.
As an example see fig.1 where $r=9,\, p=5,\,\sigma=4.$

A function $Y$ on $T$ may be represented as a vector
$Y(x)=[y_j(x)]_{j\in J}$, $x\in [0,1],$ where $J:=\{j:\;
j=\overline{1,r}\}$, and the function $y_j(x)$ is defined on the edge
$e_j$. Let $q=[q_j]_{j\in J}$ is an integrable complex-valued function
on $T$ which is called the potential.
Consider the Sturm-Liouville equation on $T$:
$$
-y''_j(x)+q_j(x)y_j(x)=\la y_j(x),\quad x\in [0,1],                       \eqno(2.1)
$$
where $j\in J$, $\la$ is the spectral parameter, the functions
$y_j(x),\, y'_j(x)$ are absolutely continuous on $[0,1]$ and satisfy
the following matching conditions in each internal vertex $v_k$,
$k=\overline{p+1,r}$:
$$
\left.\begin{array}{c}
y_k(0)=a_{kj}y_j(1)\quad \mbox{for all}\quad e_j\in R(v_k),\\[3mm]
y'_k(0)=\di\sum_{e_j\in R(v_k)}\Big(a_{kj}^1 y'_j(1)+a_{kj}^0 y_j(1)\Big),
\end{array}\right\}                                                      \eqno(2.2)
$$
where $a_{kj},\,a_{kj}^0,\, a_{kj}^1$ are complex numbers, and
$a_{kj}a_{kj}^1\ne 0$. Assume that
$$
r_k:=\di\sum_{e_j\in R(v_k)}\di\frac{a_{kj}^1}{a_{kj}}\ne -1.            \eqno(2.3)
$$
Condition (2.3) is called the regularity condition for matching.
Sturm-Liouville operators on $T$ which do not satisfy the regularity
condition for matching (2.3), possess qualitatively different properties
for the formulation and the investigation of direct inverse problems,
and are not considered here; they require a separate investigation.
We note that if $a_{kj}=a_{kj}^1=1,\,a_{kj}^0=0$ for all $k,j,$
then the conditions (2.2) are called the standard conditions.
For standard matching conditions, (2.3) is satisfied obviously.
In electrical circuits, (2.2) expresses Kirchhoff's law; in elastic
string network, it expresses the balance of tension, and so on.

Consider the linear forms in the boundary vertices $v_j,\,j\in\Gamma$:
$$
U_{js}(Y):=h_{js}^1 Y'_{|v_j}+h_{js}^0 Y_{|v_j},
\quad s=0,1,\; j=\overline{0,p},
$$
where $h_{js}^\nu$ are complex numbers such that
$\det[h_{js}^\nu]_{s,\nu=0,1}\ne 0.$ Denote by $B$ the boundary
value problem for equation (2.1) with the matching conditions
(2.2) and with the boundary conditions $U_{j0}(Y)=0$,
$j=\overline{0,p}.$ We also will consider the boundary value
problems $B_k$, $k=\overline{0,p},$ for equation (2.1) with the
matching conditions (2.2) and with the boundary conditions
$U_{k1}(Y)=0$, $U_{j0}(Y)=0$, $j=\overline{0,p}\setminus k.$

Let $\Psi_k(x,\la)=[\psi_{kj}(x,\la)]_{j\in J}$,
$k=\overline{0,p},$ be solutions of equation (2.1) satisfying
(2.2) and the boundary conditions
$$
U_{j0}(\Psi_k)=\de_{jk},\quad j=\overline{0,p},                        \eqno(2.4)
$$
where $\de_{jk}$ is the Kronecker symbol. The function $\Psi_k$
is called the Weyl solution with respect to the boundary vertex $v_k$.
Denote $M(\la)=[M_k(\la)]_{k=\overline{1,p}}$, where $M_k(\la):=
U_{k1}(\Psi_k).$ The functions $M_k(\la)$ are called the Weyl functions,
and $M(\la)$ is called the Weyl vector. We also consider the extended
Weyl vector $M^0(\la)=[M_k(\la)]_{k=\overline{0,p}}$ and the matrix
${\cal M}(\la)=[M_{jk}(\la)]_{j,k=\overline{0,p}},$ where
$M_{jk}(\la):=U_{j1}(\Psi_k).$ The vector $M^0(\la)$ is the main
diagonal of ${\cal M}(\la),$ and the vector $M(\la)$ is a part of this
diagonal.

For definiteness, we will consider the case when
$U_{j0}(Y)=Y'_{|v_j}+h_jY_{|v_j}$, $U_{j1}(Y)=Y_{|v_j}$, i.e.
$h_{j0}^1=h_{j1}^0=1,\, h_{j1}^1=0,\, h_{j0}^0=h_j$. Other cases
are treated similarly.
Let $\vv_j(x,\la),\, S_j(x,\la),\;j\in J,\, x\in[0,1]$ be solutions
of equation (2.1) on the edge $e_j$ under the initial conditions
$\vv_j(0,\la)=S'_j(0,\la)=1$, $\vv'_j(0,\la)=-h_j$, $S_j(0,\la)=0.$
For each fixed $x$, the functions $\vv_j^{(\nu)}(x,\la)$ and
$S_j^{(\nu)}(x,\la)$, $\nu=0,1,$ are entire in $\la$ of order $1/2$.
Moreover, $\langle \vv_j(x,\la), S_j(x,\la)\rangle\equiv 1,$ where
$\langle y,z\rangle :=yz'-y'z$ is the Wronskian of $y$ and $z.$
Denote $M^1_{kj}(\la)=\psi_{kj}(0,\la),$
$M^0_{kj}(\la)=\psi'_{kj}(0,\la)+h_j\psi_{kj}(0,\la).$ Then
$$
\psi_{kj}(x,\la)=M^0_{kj}(\la)S_j(x,\la)+M^1_{kj}(\la)\vv_j(x,\la).    \eqno(2.5)
$$
In particular, for $k=\overline{1,p}$, we have $M^1_{kk}(\la)=
M_{k}(\la)$, $M^0_{kk}(\la)=1$, $M^0_{kj}(\la)=0$ for
$j=\overline{1,p}\setminus k$, and consequently,
$$
\psi_{kk}(x,\la)=
S_k(x,\la)+M_{k}(\la)\vv_k(x,\la),\quad k=\overline{1,p}.             \eqno(2.6)
$$
Substituting (2.5) into (2.2) and (2.4) we obtain a linear algebraic
system $s_k$ with respect to $M^s_{kj}(\la).$ The determinant
$\Delta(\la)$ of this system is an entire function of order 1/2,
and the zeros of $\Delta(\la)$ coincide with the eigenvalues
$\Lambda=\{\la_{l}\}_{l\ge 0}$ of the boundary value problem $B$.
The function $\Delta(\la)$ is called the characteristic function
for $B.$ Solving the system $s_k$ we get by Cramer's rule:
$M_{kj}^s(\la)=\Delta_{kj}^s(\la)/\Delta(\la)$, $s=0,1,\,
j=\overline{1,r},$ where the determinant $\Delta_{kj}^s(\la)$
is obtained from $\Delta(\la)$ by the replacement of the column
which corresponds to $M_{kj}^s(\la)$ by the column of the free
terms. In particular,
$$
M_k(\la)=\Delta_k(\la)/\Delta(\la),\quad k=\overline{1,p},          \eqno(2.7)
$$
where $\Delta_k(\la):=\Delta_{kk}^0(\la).$ The function
$\Delta_k(\la)$ is entire in $\la$ of order 1/2, and its
zeros coincide with the eigenvalues $\Lambda_k=\{\la_{lk}\}_{l\ge 0}$
of the boundary value problem $B_k$. This function is called the
characteristic function for $B_k$. The eigenvalues $\la_{l}$ and
$\la_{lk}$ lie in the strip $\{\la:\; |\mbox{Im}\,\la|\le c_1,\;
\mbox{Re}\,\la \ge c_0$ for some $c_0,\,c_1$. It follows from
(2.7) that the Weyl functions $M_k(\la)$ are meromorphic in
$\la$ with the poles $\{\la_{l}\}_{l\ge 0}.$ If all poles are
simple we can introduce the data $S:=\{\la_{l},\al_{lk}\}_{l\ge 0,
\,k=\overline{1,p}}$, where $\al_{lk}$ are residues of $M_k(\la)$
at $\la_{l}$; the data $S$ are called the spectral data for $B.$
We note that it is also possible to define the spectral data in
the case of multiple poles (as the coefficients of the main part of
Laurent's series for the Weyl functions), but for simplicity here
we confine ourselves to the case of simple poles.
In particular, in the self-adjoint case (when $q$ and $h$ are real)
the Weyl functions always have only simple poles.

We study three inverse problems of recovering the potential
$q=[q_j]_{j\in J}$ and the coefficients $h=[h_j]_{j\in J}$ from
the following spectral characteristics:

{\bf Inverse problem 2.1.} Given the Weyl vector $M=[M_k]_{k=\overline{1,p}},$
construct $q=[q_j]_{j\in J}$ and $h=[h_j]_{j\in J}$.

{\bf Inverse problem 2.2.} Given the system of $p+1$ spectra: $\Lambda$
and $\Lambda_k,\,k=\overline{1,p},$ construct $q=[q_j]_{j\in J}$ and
$h=[h_j]_{j\in J}$.

{\bf Inverse problem 2.3.} From the given spectral data $S$ find
$q=[q_j]_{j\in J}$ and $h=[h_j]_{j\in J}$.

For each of these inverse problems we provide a constructive procedure
for the solution and prove its uniqueness (see section 2.3).
We note that the notion of the Weyl vector $M$ is a generalization
of the notion of the Weyl function (m-function) for the classical
Sturm-Liouville operator [12], and Inverse Problem 2.1 is a
generalization of the classical inverse problem for Sturm-Liouville
operator on an interval from the Weyl function, and (which is equivalent)
from the spectral measure. Inverse Problem 2.2 is a generalization of
Borg's inverse problem for the Sturm-Liouville operator on an
interval from two spectra. Inverse Problem 2.3 is a generalization of the
classical Marchenko's inverse problem for the Sturm-Liouville operator on
an interval (see [12] for details).

\medskip
{\bf 2.2. Properties of spectral characteristics. }
Let $\la=\rho^2,\,\mbox{Im}\,\rho\ge 0.$ Denote
$\Lambda:=\{\rho:\;\mbox{Im}\,\rho\ge 0\}$,
$\Lambda_\de:=\{\rho:\;\arg\rho\in[\de,\pi-\de]\}$, $\de>0.$
The following theorem describes asymptotical properties of the Weyl
solutions [21].

\smallskip
{\bf Theorem 2.4. }{\it Let $e_j\in{\cal E}^{(\mu)}.$ Then for $\nu=0,1,
\;\rho\in\Lambda_\de,\; |\rho|\to\iy,$ uniformly in $x\in[0,1]$, one has
$$
\psi_{0j}^{(\nu)}(x,\la)=B_j(\rho)\exp(i\rho\mu)
\Big((-i\rho)^{\nu-1}\exp(-i\rho x)[1]-
(i\rho)^{\nu-1}d_j\exp(i\rho x)[1]\Big),                           \eqno(2.8)
$$
where $[1]=1+O(\rho^{-1}),$ $d_j=1$ for $j=\overline{1,p},$
and $d_j=(1+r_j)^{-1}(1-r_j)$ for $j=\overline{p+1,r}.$
Moreover, for $\rho\in\Lambda_\de,\; |\rho|\to\iy,$ one has
$B_j(\rho)=b_j[1],\quad b_j\ne 0,\; b_{p+1}=1.$
In particular, for $\rho\in\Lambda_\de,\;|\rho|\to\iy,$}
$$
\psi_{0j}^{(\nu)}(x,\la)=(-i\rho)^{\nu-1}b_j
\exp(i\rho(\mu-x))[1],\quad \nu=0,1,\; x\in(0,1].                   \eqno(2.9)
$$

Symmetrically to (2.8)-(2.9), one can get the asymptotics for all
other Weyl solutions $\Psi_k,\,k=\overline{1,p}.$ In particular,
for $k=\overline{1,p}$, $\nu=0,1,$ one has
$$
\psi_{kk}^{(\nu)}(x,\la)=(i\rho)^{\nu-1}\exp(i\rho x)[1],
\; M_k(\la)=(i\rho)^{-1}[1],\;
\rho\in\Lambda_\de,\;|\rho|\to\iy,\; x\in[0,1).                    \eqno(2.10)
$$
Moreover,
$$
\vv_k^{(\nu)}(x,\la)=\frac{1}{2}\Big((i\rho)^{\nu}\exp(i\rho x)[1]
+(-i\rho)^{\nu}\exp(-i\rho x)[1]\Big),\quad |\rho|\to\iy.          \eqno(2.11)
$$

Fix $k=\overline{1,p},$ and consider the following auxiliary
inverse problem on the edge $e_k$, which is called IP(k).

\smallskip
{\bf IP(k).} Given $M_k(\la),$ construct $q_k(x),\, x\in[0,1]$
and $h_k$.

\smallskip
In the inverse problem IP(k) we construct the potential on the edge
$e_k$, but the Weyl function $M_k(\la)$ brings global information
from the whole graph. Therefore, IP(k) is not a local inverse problem
related only to the edge $e_k$.
Let us prove the uniqueness of the solution of the inverse problem IP(k).
For this purpose together with $T$ we consider a tree $\tilde T$ of the
same form but with different $\tilde q$ and $\tilde h.$ Everywhere below
if a symbol $\al$ denotes an object related to $T,$ then $\tilde\al$ will
denote the analogous object related to $\tilde T.$

\smallskip
{\bf Theorem 2.5. }{\it If $M_k(\la)=\tilde M_k(\la),$ then
$q_k(x)=\tilde q_k(x)$ a.e. on [0,1], and $h_k=\tilde h_k.$ Thus,
the specification of the Weyl function $M_k$ uniquely determines
the potential $q_k$ on the edge $e_k$ and the coefficient $h_k$.}

\smallskip
{\it Proof.} We introduce the functions
$$
P^k_{1s}(x,\la)=
(-1)^{s-1}\Big(\vv_k(x,\la)\tilde\psi^{(2-s)}_{kk}(x,\la)
-\tilde\vv^{(2-s)}_k(x,\la)\psi_{kk}(x,\la)\Big),\quad s=1,2.   \eqno(2.12)
$$
It follows from (2.6) that $\langle \vv_k(x,\la),
\psi_{kk}(x,\la)\rangle\equiv 1.$ Then, by direct calculations
we get
$$
\vv_k(x,\la)=P^k_{11}(x,\la)\tilde\vv_k(x,\la)
+P^k_{12}(x,\la)\tilde\vv'_k(x,\la).                            \eqno(2.13)
$$
Using (2.10)-(2.12) we obtain
$$
P^k_{1s}(x,\la)=\de_{1s}+O(\rho^{-1}),\quad
\rho\in\Lambda_\de,\;|\rho|\to\iy,\; x\in(0,1].                 \eqno(2.14)
$$
According to (2.6) and (2.12),
$$
P^k_{1s}=(-1)^{s-1}((\vv_k\tilde S^{(2-s)}_k-S_k\tilde\vv_k^{(2-s)})
+(\tilde M_k-M_k)\vv_k\tilde\vv^{(2-s)}_k).
$$
Since $M_k(\la)=\tilde M_k(\la),$ it follows that for each fixed
$x,$ the functions $P^k_{1s}(x,\la)$ are entire in $\la.$
Together with (2.14) this yields $P^k_{11}(x,\la)\equiv 1$,
$P^k_{12}(x,\la)\equiv 0.$ Substituting this relations into
(2.13) we get $\vv_k(x,\la)\equiv\tilde\vv_k(x,\la)$ for all
$x$ and $\la,$ and consequently, $q_k(x)=\tilde q_k(x)$ a.e.
on $[0,1]$, and $h_k=\tilde h_k$.
\B

\smallskip
Using the method of spectral mappings [11] for the Sturm-Liouville
operator on the edge $e_k$ one can get a constructive procedure
for the solution of the local inverse problem IP(k); for details
see [11], [21].

\smallskip
{\it Problem $Z(T,v_0,a).$} Let $\Psi=[\psi_j]_{j\in J}$ be
the solution of equation (2.1) on $T$ satisfying (2.2) and
the boundary conditions
$$
\Psi_{|v_0}=a,\quad U_{j0}(\Psi)=0,\;j=\overline{1,p},           \eqno(2.15)
$$
where $a$ is a complex number. Denote $m_j^1(\la)=
\psi_j(0,\la), \,m_j^0(\la)=\psi'_j(0,\la)+h_j\psi_j(0,\la),
\,j\in J.$ Then

\smallskip
$$
\psi_j(x,\la)=m_j^0(\la)S_j(x,\la)+m_j^1(\la)\vv_j(x,\la).       \eqno(2.16)
$$
Substituting (2.16) into (2.2) and (2.15) we obtain a linear
algebraic system with respect to $m_j^s(\la),$ $j\in J.$ The
determinant of this system is $\Delta_0(\la).$ Solving this
system we find the transition matrix
$[m_j^0(\la),m_j^1(\la)]_{j\in J}$ for $T$ with respect to
$v_0$ and $a.$ The problem of calculating the transition matrix
$[m_j^0(\la),m_j^1(\la)]_{j\in J}$ is called
Problem $Z(T,v_0,a).$ This problem will be used below for
describing the procedure for the solution of the inverse problems.

\smallskip
Let us introduce the Weyl solutions for internal vertices.
Fix $v_k\in V.$ Denote $T_k^0:=\{z\in T:\; v_k<z\}$,
$T_k:=T\setminus T_k^0$. Clearly, $T_k$ is a tree with the root $v_0.$
Let $\Gamma_k$ be the set of boundary vertices of $T_k$, and let
$E_k$ be the set of boundary edges of $T_k$. Denote $J_k:=
\{j:\; e_j\in T_k\}.$ If $Y=[y_j]_{j\in J}$ is a function on $T$,
then $\{Y\}_k:=[y_j]_{j\in J_k}$ is a function on $T_k$.

Fix $v_k\notin\Gamma$ (i.e. $k=\overline{p+1,r}$). Let
$\Psi_k(x,\la)=[\psi_{kj}(x,\la)]_{j\in J_k}$ be the solution of
equation (2.1) on $T_k$ satisfying (2.2) and the boundary conditions
$U_{j0}(\Psi_k)=\de_{kj}$, $v_j\in\Gamma_k$, where $U_{k0}(Y)=
Y'_{|v_k}+h_kY_{|v_k}$, and $h_k$ is a real number. The function
$\Psi_k$ is the Weyl solution of (2.1) on $T_k$ with respect
to the vertex $v_k$. Denote by $M_k(\la):=\psi_{kk}(0,\la)$,
$k=\overline{p+1,r}$ the Weyl functions for $T_k$ with respect
to $v_k$.

Fix $v_m\notin\Gamma.$ Let $e_k=[v_m,v_k]\in R(v_m).$ Then
$$
M_m(\la)=a_{mk}(A_{mk}(\la))^{-1}\psi_{kk}(1,\la),               \eqno(2.17)
$$
$$
A_{mk}(\la)=\di\sum_{e_j\in R(v_m)}a_{mj}^1\psi'_{kj}(1,\la)
+a_{mk}\Big(h_m+\di\sum_{e_j\in R(v_m)}
\frac{a_{mj}^0}{a_{mj}}\Big)\psi_{kk}(1,\la).                    \eqno(2.18)
$$
Denote $M^1_{kj}(\la)=\psi_{kj}(0,\la)$, $M^0_{kj}(\la)=
\psi'_{kj}(0,\la)+h_j\psi_{kj}(0,\la)$ for $k=\overline{p+1,r}$,
$j\in J_k$. Then (2.5) and (2.6) are valid for $k=\overline{1,r}$,
$j\in J_k$, where $J_k=J$ for $k=\overline{1,p}.$
In particular, this yields
$$
\psi^{(\nu)}_{kj}(1,\la)=M^0_{kj}(\la)S_j^{(\nu)}(1,\la)
+M^1_{kj}(\la)\vv_j^{(\nu)}(1,\la),\quad
\nu=0,1,\;k=\overline{1,r},\; j\in J_k,                         \eqno(2.19)
$$
$$
\psi^{(\nu)}_{kk}(1,\la)=S_k^{(\nu)}(1,\la)+
M_k(\la)\vv_k^{(\nu)}(1,\la),\quad \nu=0,1,\;k=\overline{1,r}.  \eqno(2.20)
$$

\medskip
{\bf 2.3. Solution of the inverse problems.} First consider Inverse
problem 2.1. Let the Weyl vector
$M(\la)=[M_k(\la)]_{k=\overline{1,p}}$ for the tree $T$ be given.
The procedure for the solution of Inverse Problem 1 consists in
the realization of the so-called $A_\mu$- procedures successively
for $\mu=\sigma,\sigma-1,\ldots,1,$ where $\sigma$ is the height
of the tree $T.$ Let us describe $A_\mu$- procedures.

{\bf ${\bf A_\sigma}$- procedure. } 1) For each edge
$e_k\in{\cal E}^{(\sigma)},$ we solve the local inverse problem IP(k)
and find $q_k(x),\,x\in[0,1]$ on the edge $e_k$ and $h_k$.

2) For each $e_k\in{\cal E}^{(\sigma)},$ we construct $\vv_k(x,\la),\,
S_k(x,\la),\,x\in[0,1],$ and calculate $\psi^{(\nu)}_{kk}(1,\la)$,
$\nu=0,1,$ by (2.20).

3) {\it Returning procedure.} For each fixed $v_m\in V^{(\sigma-1)}
\setminus\Gamma$ and for all $e_j,e_k\in R(v_m),\,j\ne k,$
we construct $M_{kj}^s(\la),\,s=0,1,$ by the formulae
$$
M_{kj}^0(\la)=0,\; M_{kj}^1(\la)=a_{mk}\psi_{kk}(1,\la)
(a_{mj}\vv_j(1,\la))^{-1}, \quad e_j,e_k\in R(v_m),\;j\ne k.
$$

4) For each fixed $v_m\in V^{(\sigma-1)}\setminus\Gamma$ we calculate
the Weyl function $M_m(\la)$ by (2.17), where $A_{mk}(\la)$ and
$\psi'_{kj}(1,\la)$ are constructed via (2.18) and (2.19).

\smallskip
Now we carry out $A_\mu$- procedures for $\mu=\overline{1,\sigma-1}$
by induction. Fix $\mu=\overline{1,\sigma-1},$ and suppose that
$A_\sigma,\ldots,A_{\mu+1}$- procedures have been already carried out.

\smallskip
{\bf ${\bf A_\mu}$- procedure. } For each $v_k\in V^{(\mu)},$
the Weyl functions $M_k(\la)$ are given. Indeed, if
$v_k\in V^{(\mu)}\cap\Gamma,$ then $M_k(\la)$ are given a priori,
and if $v_k\in V^{(\mu)}\setminus\Gamma,$ then $M_k(\la)$ were
calculated on the previous steps according to
$A_\sigma,\ldots,A_{\mu+1}$- procedures.

1) For each edge $e_k\in{\cal E}^{(\mu)},$ we solve the local
inverse problem IP(k) and find $q_k(x),\,x\in[0,1]$ on the edge
$e_k$ and $h_k$. If $\mu=1,$ then Inverse Problem 2.1 is solved,
and we stop our calculations. If $\mu>1,$ we go on to the next step.

2) For each $e_k\in{\cal E}^{(\mu)},$ we construct $\vv_k(x,\la),\,
S_k(x,\la),\,x\in[0,1],$ and calculate $\psi^{(\nu)}_{kk}(1,\la)$,
$\nu=0,1,$ by (2.20).

3) {\it Returning procedure.} For each fixed $v_m\in V^{(\mu-1)}
\setminus\Gamma$ and for any fixed $e_k,e_i\in R(v_m),\,i\ne k,$
we consider the tree $T^1_i:=T^0_i\cup\{e_i\}$ with the root $v_m$.
Solving the problem $Z(T^1_i,v_m,\psi_{kk}(1,\la)),$ we calculate
the transition matrix $[M_{kj}^0(\la), M_{kj}^1(\la)]$ for $e_j\in T^1_i$.

4) For each fixed $v_m\in V^{(\mu-1)}\setminus\Gamma$ we calculate
the Weyl function $M_m(\la)$ by (2.17), where $A_{mk}(\la)$ and
$\psi'_{kj}(1,\la)$ are constructed via (2.18) and (2.19).

\smallskip
Thus, we have obtained the solution of Inverse Problem 2.1 and
proved its uniqueness, i.e. the following assertion holds.

\smallskip
{\bf Theorem 2.6. }{\it The specification of the Weyl vector $M$
uniquely determines the potential $q$ on $T$ and $h.$ Solution of
Inverse Problem 2.1 can be obtained by executing successively
$A_\sigma,A_{\sigma-1},\ldots,A_{1}$ --procedures.}

\smallskip
Consider now Inverse Problem 2.2. Let the spectra $\Lambda$ and
$\Lambda_k\,$ $k=\overline{1,p}$ be given. By Hadamard's theorem,
the functions $\Delta(\la)$ and $\Delta_k(\la)$ are uniquely
determined up to multiplicative constants by their zeros:
$$
\Delta(\la)=C\prod_{l=0}^{\iy}\Big(1-\frac{\la}{\la_{l}}\Big),\quad
\Delta_k(\la)=C_k\prod_{l=0}^{\iy}\Big(1-\frac{\la}{\la_{lk}}\Big)
$$
(the case when $\Delta(0)=0$ or/and $\Delta_k(0)=0$ requires minor
modifications). Then, by virtue of (2.7),
$$
M_k(\la)=m_k\prod_{l=0}^{\iy}\Big(1-\frac{\la}{\la_{lk}}\Big)
\Big(1-\frac{\la}{\la_{l}}\Big)^{-1},\quad k=\overline{1,p},         \eqno(2.21)
$$
where $m_k$ is a constant. Using (2.10) we obtain
$$
m_k=\lim_{|\rho|\to\iy}(i\rho)^{-1}\prod_{l=0}^{\iy}
\Big(1-\frac{\la}{\la_{l}}\Big)\Big(1-\frac{\la}{\la_{lk}}\Big)^{-1},
\quad \rho\in\Lambda_\de,\; k=\overline{1,p}.                        \eqno(2.22)
$$
Thus, using the given spectra, one can construct uniquely
the Weyl vector $M(\la)=[M_k(\la)]_{k=\overline{1,p}}$ by
(2.21) and (2.22). In other words, the solution of Inverse
Problem 2.2 is reduced to the solution of Inverse Problem 2.1,
and the following assertion holds.

\smallskip
{\bf Theorem 2.7. }{\it The specification of the spectra
$\Lambda$ and $\Lambda_k$, $k=\overline{1,p}$ uniquely determines
the potential $q$ on $T$ and $h.$ For constructing the solution
of Inverse Problem 2.2 we calculate the Weyl vector $M$ by
(2.21)-(2.22), and then construct $q$ and $h$ by solving
Inverse Problem 2.1.}

\smallskip
Consider now Inverse Problem 2.3. Let all poles of $M(\la)$
are simple, and let the spectral data $S$ be given. Then it is possible
the Weyl functions via
$$
M_k(\la)=\sum_l \frac{\al_{lk}}{\la-\la_{l}},                         \eqno(2.23)
$$
where the series in (2.23) converge "with brackets" [21]. Thus, the
solution of Inverse Problem 2.3 is reduced to the solution of Inverse
Problem 2.1, and the following assertion holds.

\smallskip
{\bf Theorem 2.8. }{\it The specification of the spectral
data $S$ uniquely determines the potential $q$ on $T$ and $h.$
For constructing the solution of Inverse Problem 2.3 we
calculate the Weyl vector $M$ by (2.23), and then
construct $q$ and $h$ by solving Inverse Problem 2.1.}

\smallskip
Above mentioned results of section 2 are obtained in [14], [21].

\smallskip
{\bf Remark 2.9. } Consider the extended Weyl vector $M^0(\la)
=[M_k(\la)]_{k=\overline{0,p}}$, which include additional Weyl
function related to the root of the tree. For the self-adjoint case
the specification of $M^0(\la)$ is equivalent to the specification
of the corresponding extended spectral data
$S^0:=\{\la_{l},\al_{lk}\}_{l\ge 0,\,k=\overline{0,p}}$, where
$\al_{lk}$ are residues of the function $M_k(\la)$ at $\la_{l}$.
In [13] the uniqueness theorem of recovering the potential from the
extended spectral data $S^0$ was proved for the self-adjoint case and
for the standard matching conditions and the Dirichlet boundary conditions.
In [15] the uniqueness theorem from the matrix ${\cal M}(\la):=
[M_{jk}(\la)]_{j,k=\overline{0,p}}$ was proved.  Clearly, $M^0(\la)$
is the diagonal of the matrix ${\cal M}(\la)$, and the Weyl vector is
a part of this diagonal. Therefore, the inverse problems from $M^0(\la)$
and ${\cal M}(\la)$ are overdetermined. We also note that in [22]
a particular case of inverse problem from 1 [21] is considered,
and the uniqueness theorem was proved, which coincides with the
uniqueness theorem from [14], and also calculation of the lengths
of the edges from the matrix ${\cal M}(\la)$ is discussed.

\begin{center}
{\bf 3. Sturm-Liouville operators on compact graphs with cycles}
\end{center}

In this section we present results on inverse spectral problems for Sturm-Liouville
operators on compact graphs with cycles. The presence of cycles makes the investigation
of the inverse problems more complicated. Therefore, the solutions of inverse problems
for graphs with cycles were obtained later than for trees. In subsection 3.1 of this
review paper we consider the so-called bush-graph, i.e. an arbitrary graph having
only one cycle. Such graphs have a specific character for studying inverse problems.
In subsection 3.2 we investigate Sturm-Liouville operators on an arbitrary compact graph.
We give the statement of the inverse problem, obtain an algorithm for constructing its
solution and establish the uniqueness of the solution.

\medskip
{\bf 3.1. Inverse problems on the bush-graph.}
Consider a compact graph $G$ in ${\bf R^{\ell}}$ with the
set of edges ${\cal E}=\{e_0,\ldots, e_r\}$ and the set of vertices
$W=V\cup U,$ where $V=\{v_1,\ldots, v_r\}$, $U=\{u_1,\ldots, u_N\}.$
The graph has the form $G=e_0\cup T,$ where $e_0$ is a cycle,
$u_i\in e_0$, $i=\overline{1,N},$ $v_j\notin e_0$, $j=\overline{1,r},$
$T\cap e_0=U,$ $T=\bigcup_{j=1}^m T_j\,,$
and $T_j$ is a tree with a root from $U$ and with one rooted edge
from ${\cal E}.$ The set $T$ consists of $N$ groups of trees:
$$
T=Q_1\ldots Q_N\,,\quad Q_i\cap e_0=u_i\,,
$$
i.e. all trees from $Q_i$ have the common root $u_i$. Let $m_i$ be
the number of trees in the block $Q_i$; hence $m_1+\cdots+m_N=m.$
Denote $s_0=1,$ $s_i=m_1+\cdots+m_i$, $i=\overline{1,N}.$ Then
$$
Q_i=\bigcup_{j=s_{i-1}+1}^{s_i} T_j,\,,\quad
\bigcap_{j=s_{i-1}+1}^{s_i} T_j=u_i,\; i=\overline{1,N}.
$$
Fix $i=\overline{1,N},$ $j=\overline{1,m},$ and consider the tree $T_j\in Q_i$.
For two points $a,b\in T_j$ we will write $a\le b$ if $a$ lies on a unique
simple path connecting the root $u_i$ with $b.$ We will write $a<b$ if
$a\le b$ and $a\ne b.$
If $a<b$ we denote $[a,b]:=\{z\in T_j:\; a\le z\le b\}.$ In particular,
if $e=[v,w]$ is an edge, we call $v$ its initial point, $w$ its end
point and say that $e$ emanates from $v$ and terminates at $w.$
For each vertex $v\in T_j$ we denote by $R(v):=\{e\in T_j:\;
e=[v,w],\, w\in T_j\}$ the set of edges emanated from $v.$ For each
$v\in T_j\cap V$ we denote by $|v_j|$ the number of edges between $u_i$
and $v.$ The number $|v|$ is called the order of $v.$
The cycle $e_0$ consists of $N$ parts:
$$
e_0=e_1^0\ldots e_N^0\,,\quad e_i^0=[u_i,u_{i+1}],\;
i=\overline{1,N},\; u_{N+1}:=u_1.
$$
For definiteness we enumerate the vertices $v_j\in V$ as follows:
$\Gamma:=\{v_1,\ldots, v_p\}$ are boundary vertices of $G$,
and $v_j,\, j>p$ are enumerated in order of increasing $|v_j|$.
We enumerate the edges similarly, namely:
$e_j=[v_{j_k},v_j]$, $j=\overline{1,r},\, j_k<j.$ In particular,
$E:=\{e_1,\ldots, e_{p}\}$ is the set of boundary edges of $G$.

\begin{center}
\unitlength=1mm
\begin{picture}(150,80)

\qbezier(35,35)(36.15,48.85)(50,50)

\qbezier(50,50)(63.85,48.85)(65,35)

\qbezier(65,35)(63.85,21.15)(50,20)

\qbezier(50,20)(36.15,21.15)(35,35)

\put(39.25,45.75){\line(-1,0){22}} \put(39.25,45.75){\circle*{0.7}}

\put(39.25,45.75){\line(-1,1){15}}                                   \put(24.25,60.75){\circle*{0.7}}

\put(17.25,45.75){\line(-3,2){15}} \put(17.25,45.75){\circle*{0.7}}  \put(02.25,55.75){\circle*{0.7}}

\put(17.25,45.75){\line(-3,-2){17}}                                  \put(00.25,34.45){\circle*{0.7}}

\put(39.25,24.25){\line(-1,-1){20}} \put(39.25,24.25){\circle*{0.7}} \put(19.25,04.25){\circle*{0.7}}

\put(64.60,38){\line(1,-2){14}}    \put(64.60,38){\circle*{0.7}}     \put(78.60,10){\circle*{0.7}}

\put(64.60,38){\line(3,-1){27}}

\put(64.60,38){\line(3,2){22}}

\put(86.60,52.62){\line(2,5){9}} \put(86.60,52.62){\circle*{0.7}}    \put(95.60,75.00){\circle*{0.7}}

\put(86.60,52.62){\line(5,1){27}}                                    \put(113.60,58){\circle*{0.7}}

\put(91.60,28.85){\line(2,-3){17}} \put(91.60,28.85){\circle*{0.7}}  \put(108.60,03.30){\circle*{0.7}}

\put(91.60,28.85){\line(3,2){20}}

\put(111.60,42.2){\line(1,1){20}} \put(111.60,42.2){\circle*{0.7}}   \put(131.60,62.2){\circle*{0.7}}

\put(111.60,42.2){\line(6,-1){35}}                                   \put(146.60,36.4){\circle*{0.7}}

\put(111.60,42.2){\line(2,-3){16}}                                   \put(127.60,18.25){\circle*{0.7}}

\put(33,24){\small $u_1$}

\put(23.5,14){\small $e_1$}

\put(14.0,02.5){\small $v_1$}

\put(30,35){\small $e^0_1$}

\put(50,51.5){\small $e^0_2$}

\put(71,47.5){\small $e_{13}$}

\put(77,35){\small $e_{14}$}

\put(78,7){\small $v_{11}$}

\put(100,16){\small $e_{10}$}

\put(109,1){\small $v_{10}$}

\put(131.5,63.5){\small $v_7$}

\put(130,40.2){\small $e_8$}

\put(147.3,35){\small $v_8$}

\put(97,57){\small $e_6$}

\put(113,59.5){\small $v_6$}

\put(87,66){\small $e_5$}

\put(95.5,75.5){\small $v_5$}

\put(16,42){\small $v_{12}$}

\put(9,51.5){\small $e_3$}

\put(-2,57.5){\small $v_3$}

\put(4,41){\small $e_2$}

\put(24,47){\small $e_{12}$}

\put(36.5,48.5){\small $u_2$}

\put(31,54){\small $e_4$}

\put(20,62.5){\small $v_4$}

\put(-4.5,32.5){\small $v_2$}

\put(54,24){\small $e^0_3$}

\put(63.5,41.5){\small $u_3$}

\put(72,23){\small $e_{11}$}

\put(88,31.5){\small $v_{14}$}

\put(120,29){\small $e_9$}

\put(115.7,52){\small $e_7$}

\put(106,44){\small $v_{15}$}

\put(100,33){\small $e_{15}$}

\put(80,54){\small $v_{13}$}

\put(128,16){\small $v_9$}

\put(70,0){\small \llap{fig. 3.1}}

\end{picture}

\end{center}

Let $d_j$ be the length of the edge $e_j$, $j=\overline{0,r}.$
Each edge $e_j$, $j=\overline{0,r}$ is viewed as a segment $[0,d_j]$
and is parameterized by the parameter $x_j\in [0,d_j].$ It is convenient
for us to choose the following orientation: for $j=\overline{1,r}$
the end vertex $v_j$ corresponds to $x_j=0,$ and the initial vertex
corresponds to $x_j=d_j$; for the cycle $e_0$ both ends $x_0=+0$
and $x_0=d_0-0$ correspond to the point $u_1$. Let $d_i^0$ be
the length of $e_i^0$. Then $d_0=d_1^0+\cdots+d_N^0$. Each part
$e_i^0$ ($i=\overline{1,N}$) of $e_0$ is parameterized by the
parameter $\xi_i\in [0,d_i^0],$ where $\xi_i=0$ corresponds
to the point $u_i$, and $\xi_i=d_i^0$ corresponds to $u_{i+1}$.
As an example see figure 3.1 where $N=3,\,r=15,\, p=11,\,\sigma=3,\,
m=6,\,m_1=1,\, m_2=2,\, m_3=3.$

A function $Y$ on $G$ may be represented as
$Y=\{y_j\}_{j=\overline{0,r}}$, where the function $y_j(x_j),$
$x_j\in [0,d_j],$ is defined on the edge $e_j$.
The function $y_0$ may be represented as $y_0=\{y_i^0\}_{i=
\overline{1,N}}$, where the function $y_i^0(\xi_i),$
$\xi_i\in [0,d_i^0],$ is defined on $e_i^0$.
Let $q=\{q_j\}_{j=\overline{0,r}}$ be an integrable real-valued
function on $G$; $q$ is called the potential.
Consider the differential equation on $G$:
$$
-y''_j(x_j)+q_j(x_j)y_j(x_j)=\la y_j(x_j),\quad x_j\in [0,d_j],   \eqno(3.1)
$$
where $j=\overline{0,r},$ $\la$ is the spectral parameter,
the functions $y_j, y'_j,$ $j=\overline{0,r},$ are absolutely
continuous on $[0,d_j]$ and satisfy the following matching
conditions in the internal vertices $u_i$, $i=\overline{1,N}$
and $v_k$, $k=\overline{p+1,r}$:
For $k=\overline{p+1,r}$,
$$
y_j(d_j)=y_k(0)\; \mbox{for all}\; e_j\in R(v_k),
\quad \di\sum_{e_j\in R(v_k)} y'_j(d_j)=y'_k(0),                  \eqno(3.2)
$$
and for $i=\overline{1,N}$,
$$
y_i^0(0)=y_{i-1}^0(d_{i-1}^0)=(Y_j)_{|u_i}\;
\mbox{for all}\; T_j\in Q_i,\qquad
(y_i^0)'(0)=(y_{i-1}^0)'(d_{i-1}^0)+
\di\sum_{T_j\in Q_i} (Y_j)'_{|u_i},                              \eqno(3.3)
$$
where $y_0^0:=y_N^0,\; d_0^0:=d_N^0$, $Y_j:=\{Y\}_{T_j}$.
Matching conditions (3.2)-(3.3) are called the standard conditions.
Let us consider the boundary value problem $L_0(G)$ for equation
(3.1) with the matching conditions (3.2)-(3.3) and with the
Dirichlet boundary conditions at the boundary vertices
$v_1,\ldots, v_p$:
$$
y_j(0)=0,\quad j=\overline{1,p}.
$$
Denote by $\Lambda_0=\{\la_{n0}\}_{n\ge 1}$ the eigenvalues
(counting with multiplicities) of $L_0(G).$ We also consider
the boundary value problems $L_{\nu_1,\ldots,\nu_\ga}(G),$
$\ga=\overline{1,p},$ $1\le\nu_1<\ldots\nu_\ga\le p$ for
equation (3.1) with the matching conditions (3.2)-(3.3) and
with the boundary conditions
$$
y'_i(0)=0,\; i=\nu_1,\ldots,\nu_\ga, \quad
y_j(0)=0,\; j=\overline{1,p},\; j\ne\nu_1,\ldots,\nu_\ga.
$$
Denote by $\Lambda_{\nu_1,\ldots,\nu_\ga}:=
\{\la_{n,\nu_1,\ldots,\nu_\ga}\}_{n\ge 1}$ the eigenvalues
(counting with multiplicities) of the boun\-dary value problem
$L_{\nu_1,\ldots,\nu_\ga}(G)$.
Let $S_j(x_j,\la),\; C_j(x_j,\la),\; j=\overline{0,r},\;
x_j\in[0,d_j],$ be the solutions of equation (3.1) on the edge
$e_j$ with the initial conditions
$S_j(0,\la)=C'_j(0,\la)=0,$ $S'_j(0,\la)=C_j(0,\la)=1.$
Denote $h(\la):=S_0(d_0,\la),$ $H(\la):=C_0(d_0,\la)-S'_0(d_0,\la).$
Let $\{z_n\}_{n\ge 1}$ be zeros of the entire function $h(\la),$
and put $\om_n:=\mbox{sign}\,H(z_n),$ $\Omega=\{\om_n\}_{n\ge 1}$.

We choose and fix one boundary vertex $v_{\xi_i}\in Q_i$ from each
block $Q_i$, $i=\overline{1,N},$ i.e. $s_{i-1}+1\le\xi_i\le s_i$.
Denote by $\xi:=\{k:\; k=\xi_1,\ldots,\xi_N\}$ the set of indices
$\xi_i$, $i=\overline{1,N}.$ The inverse problem is formulated as follows.

\smallskip
{\bf Inverse problem 3.1.} Given $2^N+p-N$ spectra $\Lambda_j,\;
j=\overline{0,p}$, $\Lambda_{{\nu_1,\ldots,\nu_\ga}}$,
$\ga=\overline{2,N},$ $1\le\nu_1<\ldots<\nu_\ga\le p,$ $\nu_j\in\xi,$
and $\Omega,$ construct the potential $q$ on $G.$

\smallskip
{\bf Theorem 3.2. }{\it The specification of $\Lambda_j,$
$j=\overline{0,p}$, $\Lambda_{{\nu_1,\ldots,\nu_\ga}}$,
$\ga=\overline{2,N},$ $1\le\nu_1<\ldots<\nu_\ga\le p,$
$\nu_j\in\xi,$ and $\Omega,$ uniquely determines the potential
$q$ on $G.$}

\smallskip
This uniqueness theorem is proved in [16], where an algorithm for the
solution of inverse problem 3.1 is also provided.

\smallskip
{\it Examples. } 1) Let $N=1.$ Then we specify $p+1$ spectra
$\Lambda_j,\; j=\overline{0,p},$ and $\Omega.$\\
2) Let $\sigma=1,\;N=r.$ Then we specify $2^N$ spectra
$\Lambda_0$, $\Lambda_{\nu_1,\ldots,\nu_\ga}$,
$\ga=\overline{1,N},$ $1\le\nu_1<\ldots\nu_\ga\le N$, and $\Omega.$\\
3) Consider the graph on figure 3.2. Then $N=2,\; p=4,\;
r=5,\; \sigma=2,\; m=3,\; m_1=1,\; m_2=2.$ Fix $j=2\vee 3\vee 4$.
Then we specify $\Omega$ and the following 6 spectra: $\Lambda_{0},
\Lambda_{1}, \Lambda_{2}, \Lambda_{3}, \Lambda_{4},\Lambda_{1j}$.\\
4) Consider the graph on figure 3.1. Then $N=3,\; p=11,\;r=15,\;
\sigma=3,\; m=6,\; m_1=1,\; m_2=2,\;m_3=3.$ Take, for example, $\xi_1=1,\;
\xi_2=3,\;\xi_3=6.$ Then we specify $\Omega$ and the following spectra:
$\Lambda_{13}, \Lambda_{16}, \Lambda_{36}, \Lambda_{j}, j=\overline{0,11}$.

\begin{center}
\unitlength=1mm
\begin{picture}(150,68)

\qbezier(50,35)(51.15,48.85)(65,50)

\qbezier(65,50)(78.85,48.85)(80,35)

\qbezier(80,35)(78.85,21.15)(65,20)

\qbezier(65,20)(51.15,21.15)(50,35)

\put(50,35){\line(-1,0){43}} \put(50,35){\circle*{0.7}} \put(7,35){\circle*{0.7}}

\put(79.60,38){\line(3,-1){27}}  \put(79.60,38){\circle*{0.7}}  \put(104.60,63){\circle*{0.7}}

\put(79.60,38){\line(1,1){25}}

\put(106.60,28.85){\line(1,-1){22}} \put(106.60,28.85){\circle*{0.7}}  \put(128.5,7){\circle*{0.7}}

\put(106.60,28.85){\line(3,2){28}}   \put(134.60,47.5){\circle*{0.7}}

\put(44.5,37){\small $u_2$}

\put(63,52){\small $e^0_1$}

\put(88,52){\small $e_2$}

\put(92,35){\small $e_5$}

\put(117.5,18,5){\small $e_4$}

\put(129,5){\small $v_4$}

\put(1.5,34){\small $v_1$}

\put(25,37){\small $e_1$}

\put(61.5,22.5){\small $e^0_2$}

\put(78,42.5){\small $u_1$}

\put(103.5,31.5){\small $v_5$}

\put(135,48){\small $v_3$}

\put(117,40){\small $e_3$}

\put(104,64){\small $v_2$}

\put(70,0){\small \llap{fig. 3.2}}

\end{picture}

\end{center}

{\bf 3.2. Inverse problems on an arbitrary compact graph.}
Consider a compact connected graph $G$ in ${\bf R^{\ell}}$ with the set
of edges ${\cal E}=\{e_1,\ldots, e_s\}$, the set of vertices
${\cal V}=\{v_1,\ldots, v_m\}$, and with the map $\sigma$ which
assings to each edge $e_j\in{\cal E}$ an ordered pair of (possibly equal)
vertices: $\sigma(e_j):=[u_{2j-1}, u_{2j}],$ $u_j\in{\cal V}.$
The vertices $u_{2j-1}=:\sigma^{-}(e_j)$ and $u_{2j}=:\sigma^{+}(e_j)$ are
called the {\it initial} and {\it terminal} vertices of $e_j$, respectively.
We will say that the edge $e_j$ {\it emanates} from $u_{2j-1}$ and
{\it terminates} at $u_{2j}$. The points $U:=\{u_j\}_{j=\overline{1,2s}}$
are called the {\it endpoints} for ${\cal E}.$ Each vertex $v\in{\cal V}$
generates the equivalence class (which is denoted by the same symbol $v$):
$v=\{u_{j_1},\ldots,u_{j_\nu}\}$ such that $v=u_{j_1}=\ldots=u_{j_\nu}$.
In other words, the set of endpoints $U$ is divided into $m$ equivalence
classes $v_1,\ldots,v_m$. The number of endpoints in the class $v_k$
is called the {\it valency} of $v_k$, and is denoted by $val\,(v_k)$.
The vertex $v_k\in{\cal V}$ is called a {\it boundary vertex} if
$val\,(v_k)=1.$ All other vertices are called {\it internal}.
Let ${\cal V}_0=\{v_1,\ldots, v_p\}$ be the boundary vertices, and
let ${\cal V}_1=\{v_{p+1},\ldots, v_m\}$ be the internal vertices.
The edge $e_j$ is called a {\it boundary edge} if one of its endpoints
belongs to ${\cal V}_0$. All other edges are called internal. Let
${\cal E}_0=\{e_1,\ldots, e_p\}$ be the boundary edges, and $v_k\in e_k$
for $k=\overline{1,p}$. An edge $e_k\in {\cal E}$ is called {\it adjacent}
to $v\in {\cal V},$ if $v\in e_k$. Denote by $R(v,G)$ the set of edges of
$G$ which are adjacent to $v.$ Let $l_j$ be the length of the edge $e_j$.
Each edge $e_j\in{\cal E}$ is parameterized by the parameter $x_j\in[0,l_j]$
such that the initial point $u_{2j-1}$ corresponds to $x_j=0,$ and the
terminal point $u_{2j}$ corresponds to $x_j=l_j$.

\begin{center} \unitlength=1mm
\begin{picture}(150,100)

\qbezier(35,65)(36.15,78.85)(50,80)
\qbezier(50,80)(63.85,78.85)(65,65)
\qbezier(65,65)(63.85,51.15)(50,50)
\qbezier(50,50)(36.15,51.15)(35,65)

\qbezier(65,65)(66.15,83.85)(85,85)
\qbezier(85,85)(103.85,83.85)(105,65)
\qbezier(105,65)(103.85,46.15)(85,45)
\qbezier(85,45)(66.15,46.15)(65,65)

\put(35,65){\line(-3,2){30}}
\put(05,85){\circle*{0.7}}

\put(35,65){\line(-3,-2){30}}
\put(05,45){\circle*{0.7}}

\put(35,65){\circle*{0.7}}

\put(65,65){\line(1,-6){8}}
\put(65,65){\circle*{0.7}}
\put(73,17){\circle*{0.7}}

\put(103,75){\line(2,-1){34}}
\put(103,75){\circle*{0.7}}
\put(137,58){\circle*{0.7}}

\put(103,75){\line(-1,-1){28}}
\put(75.25,47.25){\circle*{0.7}}

\put(71,14){\small $v_1$}
\put(01,42){\small $v_2$}
\put(02,87){\small $v_3$}
\put(137,55){\small $v_4$}
\put(30,68){\small $v_5$}
\put(62,74){\small $v_6$}
\put(103,76){\small $v_7$}
\put(73,44){\small $v_8$}

\put(64,35){\small $e_1$}
\put(20,53){\small $e_2$}
\put(20,75){\small $e_3$}
\put(120,67){\small $e_4$}
\put(48,81){\small $e_5$}
\put(48,52){\small $e_6$}
\put(85,62){\small $e_7$}
\put(83,86){\small $e_8$}
\put(67,55){\small $e_9$}
\put(83,47){\small $e_{10}$}

\put(70,0){\small \llap{fig. 3.3}}

\end{picture}
\end{center}

A chain of edges $\{e_{\nu_1},\ldots,e_{k,\nu_\eta}\}$ is called a
{\it cycle} if it forms a closed curve. The edge $e_j\in{\cal E}$ is
called a {\it simple} edge, if it is not a part of cycles. In particular,
all boundary edges $e_1,\ldots,e_p$ are simple.
We enumerate the edges as follows: ${\cal E}_1=\{e_1,\ldots, e_r\}$ are
simple edges, ${\cal E}_2=\{e_{r+1},\ldots, e_s\}$ are the edges which
form the set of cycles.
For definiteness, let $p>1$ (the cases $p=0$ and $p=1$ require slight
modifications; see the remark at the end of section 3.2). Let us take the
boundary vertex $v_p$ as a root. The corresponding edge $e_p$ is called the
{\it rooted edge}. For definiteness, we agree that if $e_j\in{\cal E}_1$ is
a simple edge, then the endpoint $u_{2j}$ is nearer to the root than $u_{2j-1}$.
If we contract each cycle to a point, then we get a new graph $G^{*}$ with
the set of edges ${\cal E}_1$. Clearly, $G^{*}$ is a tree (i.e. a graph
without cycles). Fix $e_k\in G^{*}.$ The minimal number $\om_k$ of edges
of $G^{*}$ between the rooted edge and $e_k$ (including $e_k$) is called
the {\it order} of $e_k$. The order of the rooted edge is equal to zero.
The number $\om:=\di\max_{e_k\in G^{*}} \om_k$ is called the order of
$G^{*}.$ We denote by ${\cal E}^{(\mu)}$, $\mu=\overline{0,\om}$, the
set of simple edges of order $\mu.$

For example, for the graph on fig. 3.3 we have
$s=10,\,m=8,\,r=p=4,\,\om=1,$ $v_1=\{u_1\},\; v_2=\{u_3\},\; v_3=\{u_5\},
\; v_4=\{u_8\},\; v_5=\{u_4,u_6,u_{10},u_{11}\},\; v_6=\{u_2, u_{9},
u_{12}, u_{16}, u_{17}\},\; v_7=\{u_7, u_{14}, u_{15}, u_{20}\},\;
v_8=\{u_{13}, u_{18}, u_{19}\}$.

A function $Y$ on $G$ may be represented as $Y=\{y_j\}_{j=\overline{1,s}}$,
where the function $y_j(x_j),$ $x_j\in [0,l_j],$ is defined on the edge
$e_j$. Denote $Y_{|u_{2j-1}}:=y_j(0),\;Y_{|u_{2j}}:=y_j(l_j),\;
\partial Y_{|u_{2j-1}}:=y'_j(0),\;\partial Y_{|u_{2j}}:=-y_j'(l_j).$
If $v\in{\cal V},$ then $Y_{|v}=0$ means $Y_{|u_j}=0$ for all $u_j\in v.$
Let $q=\{q_j\}_{j=\overline{1,s}}$ be an integrable real-valued
function on $G$; $q$ is called the potential.
Consider the following differential equation on $G$:
$$
-y''_j(x_j)+q_j(x_j)y_j(x_j)=\la y_j(x_j),\quad x_j\in [0,l_j],     \eqno(3.4)
$$
where $j=\overline{1,s},$ $\la$ is the spectral parameter,
the functions $y_j, y'_j,$ $j=\overline{1,s},$ are absolutely
continuous on $[0,l_j]$ and satisfy the following matching
conditions (MC) in each internal vertex $v_\xi\in{\cal V}_1$:
$$
Y_{|u_i}=Y_{|u_j}\; \mbox{for all}\; u_i,\,u_j\in v_\xi,
\qquad \di\sum_{u_i\in v_\xi} \partial Y_{|u_i}=0.                  \eqno(3.5)
$$
Fix $e_k\in {\cal E}_2$ and $\ee_k=0 \vee 1.$ Denote
$w_k:=u_{2k-\ee_k}$. If (3.5) holds for the set
$U\setminus\{w_k\}$, we will call these conditions the $w_k$-MC.

Let us consider the boundary value problem $L_0(G)$ for equation
(3.4) with MC (3.5) at the internal vertices ${\cal V}_1$ and with the
Dirichlet boundary conditions at the boundary vertices ${\cal V}_0$:
$$
Y_{|v_j}=0,\quad j=\overline{1,p}.                                 \eqno(3.6)
$$
We also consider the boundary value problems $L_k(G),$
$k=\overline{1,p-1},$ for equation (3.4) with MC (3.5) and
with the boundary conditions
$\partial Y_{|v_k}=0,\; Y_{|v_j}=0,\; j=\overline{1,p}\setminus k.$
Let $\Lambda_k=\{\la_{kn}\}_{n\ge 1}$, $k=\overline{0,p-1},$ be
the eigenvalues (counting with multiplicities) of $L_k(G).$

Let $L_\nu^{\xi}(G),\; \xi=\overline{r+1,s},\;\nu=0,1,$ be the
boundary value problem for equation (3.4) with $w_{\xi}$-MC and
with the boundary conditions
$\partial^\nu Y_{|w_\xi}=0,\; Y_{|v_j}=0,\; j=\overline{1,p},$
where $\partial^0 Y:=Y,\; \partial^1 Y:=\partial Y.$
We denote by $\Lambda_\nu^\xi=\{\la_{\nu n}^{\xi}\}_{n\ge 1}$ the
eigenvalues (counting with multiplicities) of $L_\nu^\xi(G).$
The inverse problem is formulated as follows.

\medskip
{\bf Inverse problem 3.3. } Given the spectra $\Lambda_k,\;
k=\overline{0,p-1},$ and $\Lambda_{\nu}^\xi,\; \xi=\overline{r+1,s},\;
\nu=0,1,$ construct the potential $q$ on $G.$

For inverse problem 3.3 the uniqueness theorem is proved (Theorem 3.4),
and a procedure for constructing the solution is obtained (Algorithm 3.6).

{\bf Theorem 3.4. }{\it The specification of the spectra
$\Lambda_k\; k=\overline{0,p-1},$ $\Lambda_{\nu}^\xi\;
\xi=\overline{r+1,s},\; \nu=0,1,$ uniquely determines the potential
$q$ on $G.$}

\smallskip
Let $S_j(x_j,\la),\; C_j(x_j,\la),\; j=\overline{1,s},\;x_j\in[0,l_j],$
be the solutions of equation (3.4) on the edge $e_j$ with the initial
conditions $S_j(0,\la)=C'_j(0,\la)=0,\; S'_j(0,\la)=C_j(0,\la)=1.$
Let $Y=\{y_j\}_{j=\overline{1,s}}$ be a solution of equation (3.4)
on $G.$ Then
$$
y_j(x_j,\la)=a_{j1}(\la)C_j(x_j,\la)+a_{j2}(\la)S_j(x_j,\la),
\quad j=\overline{1,s},                                            \eqno(3.7)
$$
where $a_{j1}(\la)$ and $a_{j2}(\la)$ do not depend on $x_j$.
Substituting (3.7) into (3.5) and (3.6), we get a linear algebraic system
$s_0$ with respect to $a_{j1}(\la),\,a_{j2}(\la),\; j=\overline{1,s}.$
The determinant $\Delta_0(\la, G)$ of system $s_0$ is an entire
function of order $1/2$. The zeros of $\Delta_0(\la, G)$ coincide
with the eigenvalues of $L_0(G).$ The function $\Delta_0(\la, G)$
is called the {\it characteristic function} for $L_0(G).$
Analogously we define the characteristic functions $\Delta_k(\la, G),$
$k=\overline{1,p-1}$ of the boundary value problems $L_k(G).$
The function $\Delta_k(\la, G)$ is obtained from $\Delta_0(\la, G)$
by the replacement of $S_k^{(\nu)}(l_k,\la),\;\nu=0,1,$ with
$C_k^{(\nu)}(l_k,\la).$ We denote by $\Delta_\nu^\xi(\la, G)$
the characteristic function of $L_\nu^\xi(G).$

Let $e_k,\; k=\overline{p,r},$ be a fixed simple edge of $G$, and let
$v=\sigma^{-}(e_k)\in{\cal V}$ be the initial point of $e_k$. The vertex
$v$ divides the graph $G$ into two parts: $G=Q\cup\hat G,$ where
$Q\cap\hat G=v,\; \hat G\cap R(v,G)=e_k$.
Consider the boundary value problem $L_0(Q,v)$ for equation (3.4) on $Q$
with MC (3.5) for $v_\xi\in{\cal V}_1\setminus\{v\}$ and with the boundary
conditions $Y_{|v_j}=0,\; v_j\in{\cal V}_0\cup\{v\}.$ Let $\Delta(\la,Q,v)$
be the characteristic function for $L_0(Q,v).$
Expanding the determinant $\Delta_0(\la, G)$ of system $s_0$ with respect
to the columns corresponding to $a_{j1}(\la)$ and $a_{j2}(\la)$ for
$e_j\in\hat G,$, we obtain the following relations.

{\it 1) } Let $v$ be an internal vertex of $Q.$ Then
$$
\Delta_0(\la, G)=\Delta_0(\la,\hat G)\Delta_0(\la, Q)
+\Delta_k(\la, \hat G)\Delta_0(\la, Q, v),                         \eqno(3.8)
$$
$$
\Delta_j(\la, G)=\Delta_0(\la, \hat G)\Delta_j(\la, Q)
+\Delta_k(\la,\hat G)\Delta_j(\la, Q, v),
\quad e_j\in{\cal E}_0\cap Q.                                      \eqno(3.9)
$$

{\it 2) } Let $v$ be a boundary vertex of $Q,$ i.e.
$R(v,Q)=: e_i\in{\cal E}_1$ consists of one simple edqe $e_i$,
and $v=\sigma^{+}(e_i).$ Then
$$
\Delta_0(\la, G)=\Delta_0(\la, \hat G)\Delta_i(\la, Q)
+\Delta_k(\la, \hat G)\Delta_0(\la, Q),                            \eqno(3.10)
$$
$$
\Delta_j(\la, G)=\Delta_0(\la, \hat G)\Delta_{ij}(\la, Q)
+\Delta_k(\la, \hat G)\Delta_j(\la, Q),\quad e_j\in{\cal E}_0\cap Q. \eqno(3.11)
$$
Here $\Delta_{ij}(\la, Q)$ is the characteristic function of the
boundary value problem $L_{ij}(Q)$ which is obtained from $L_0(Q)$ by the
replacement of the Dirichlet boundary conditions at the boundary vertices
$\sigma^{-}(e_j)$ and $\sigma^{+}(e_i),$ with the Neumann conditions.
Therefore, $\Delta_{ij}(\la, Q)$ is obtained from $\Delta_{0}(\la, Q)$
by the replacement of $S_j^{(\nu)}(l_j,\la),\; \nu=0,1,$ with
$C_j^{(\nu)}(l_j,\la)$ and by the replace\-ment of $S_i(l_i,\la),\;
C_i(l_i,\la),$ with $S'_i(l_i,\la),\; C'_i(l_i,\la),$ respectively.

Let $\la=\rho^2,\;\mbox{Im}\,\rho\ge 0.$ Denote
$\Lambda:=\{\rho:\;\mbox{Im}\,\rho\ge 0\},$ $\Lambda^\de:=
\{\rho:\;\mbox{arg}\,\rho\in[\de,\pi-\de]\}$. Let
$\lambda^0_{kn}=(\rho^0_{kn})^2,\; n\ge 1$ be the eigenvalues of the
boundary value problem $L_k(G)$ with the zero potential $q=0$.
This problem will be denoted by $L_k^0(G.)$ Let $\Delta_k^0(\la,G)$
be the characteristic function for $L_k^0(G).$ Then

1) There exists $h>0$ such that the numbers $\lambda_{kn}=\rho_{kn}^2$
lie in the strip $|\mbox{Im}\,\rho|<h.$

2) For $\rho\in\Lambda^\de,\; |\rho|\to\iy,$ one has
$\Delta_k(\la, G)=\Delta_k^0(\la, G)(1+O(\rho^{-1})).$

3) For $n\to\infty,$
$\rho_{kn}=\rho_{kn}^0 + O((\rho_{kn}^0)^{-1}).$

The characteristic functions $\Delta^{\xi}_{\nu}(\la, G)$
have similar properties. Denote
$$
\mu_{kn}^{0}=
\left\{ \begin{array}{ll}
\la_{kn}^{0},\quad & \mbox{если}\quad \la_{kn}^{0}\ne 0,\\[0.2cm]
1,\quad & \mbox{если}\quad \la_{kn}^{0}=0.
\end{array}\right.,\qquad
\mu_{kn}=
\left\{ \begin{array}{ll}
\la_{kn},\quad & \mbox{если}\quad \la_{kn}\ne 0,\\[0.2cm]
1,\quad & \mbox{если}\quad \la_{kn}=0.
\end{array}\right.
$$
Using Hadamard's factorization theorem and the properties of the
characteristic functions, one can show that the specification of
the spectrum $\Lambda_k=\{\la_{kn}\}_{n\ge 1}$ uniquely determines
the charac\-teristic function $\Delta_k(\la, G)$ by the formula
$$
\Delta_k(\la, G)=
A_k^0 \prod_{n=1}^{\iy} \frac{\la_{kn}-\la}{\mu_{kn}^{0}}\,,
\quad A_k^0=(-1)^{s_k}\frac{1}{s_k!}
\Big(\frac{\partial^{s_k}}{\partial\la^{s_k}}
\,\Delta_k^0(\la, G)\Big)_{|\la=0}\,,
$$
where $s_k\ge 0$ is the multiplicity of the zero eigenvalue of $L_k^0(G)$.
Analogously, the specification of the spectrum $\Lambda_\nu^{\xi}=
\{\la_{\nu n}^{\xi}\}_{n\ge 1}$ uniquely determines the characteristic
function $\Delta_\nu^{\xi}(\la, G).$

\smallskip
Fix $k=\overline{1,p-1},$ and consider the following auxiliary inverse
problem on the edge $e_k$, which is called $IP(e_k, G).$

$IP(e_k, G).$ {\it Даны $\Delta_0(\la, G)$ и $\Delta_k(\la, G),$
построить потенциал $q$ на $e_k$.}

{\bf Theorem 3.5. }{\it Fix $k=\overline{1,p-1}.$ The specification of
the functions $\Delta_0(\la, G)$ and $\Delta_k(\la, G)$
uniquely determines the potential $q_k$ on the edge $e_k$.}

This theorem is proved in [20], where a constructive procedure for the
 solution of the inverse problem $IP(e_k, G)$ is also obtained.

\smallskip
Fix $\xi=\overline{r+1,s}.$ Consider the following auxiliary inverse
problem, which is called $IP(e_\xi, G).$

\smallskip
$IP(e_\xi, G).$ Given $\Delta_0^\xi(\la, G)$ and
$\Delta_1^\xi(\la, G),$ construct the potential $q$ on $e_{\xi}$.

If we move slightly the endpoint $w_\xi$ to a point $w_\xi^0\notin G$
(without moving the other endpoints of $G$ and without changing the
length of $e_\xi$), then instead of $e_\xi$ we get the graph $G^{\xi}$
with the edge $e_\xi^0$ instead of $e_\xi$. The edge $e_\xi^{0}$ is a
boundary edge for $G^{\xi}$, and $w_\xi^0$ is a boundary vertex of
$G^{\xi}$. Then the inverse problem $IP(e_\xi, G)$ is equivalent to
the inverse problem $IP(e_{\xi}^0, G^{\xi}).$ Therefore, $IP(e_\xi, G),$
$\xi=\overline{r+1,s}$ is solved by the same arguments as $IP(e_k, G).$

\smallskip
{\bf Descent procedure. } Fix $k=\overline{p,r}.$ Let $e_k\in{\cal E}^{(\mu)}$
be a fixed simple edge of order $\mu$, and let $v=\sigma^{-}(e_k)\in{\cal V}$
be the initial point of $e_k$. The vertex $v$ divides the graph $G$ into two
parts $G=Q\cup\hat G,$ where $Q\cap\hat G=v,\; \hat G\cap R(v,G)=e_k$.
Then (3.8)-(3.9) or (3.10)-(3.11) hold.
Suppose that the potential $q$ is known a priori on $Q.$ Fix $e_j\in
{\cal E}_0\cap Q.$ Let $\Delta_0(\la,G)$ and $\Delta_j(\la,G)$ be given.\\
1) Solving the algebraic system (3.8)-(3.9) or (3.10)-(3.11) we calculate
$\Delta_0(\la, \hat G)$ and $\Delta_k(\la, \hat G).$\\
2) Solving the inverse problem $IP(e_k, \hat G),$ we construct
the potential $q$ on $e_k$.\\
This procedure of calculating $q$ on $e_k$ is called the
{\it descent procedure}.

\smallskip
{\it Solution of Inverse problem 3.3. }
Let the spectra $\Lambda_k,\; k=\overline{0,p-1}$ and
$\Lambda_{\nu}^\xi,\; \xi=\overline{r+1,s},\; \nu=0,1,$ be given.
Then one can uniquely construct the characteristic functions
$\Delta_k(\la, G),\; k=\overline{0,p-1}$ and $\Delta_\nu^\xi(\la, G),
\; \xi=\overline{r+1,s},\; \nu=0,1.$
The solution of Inverse problem 1 can be found by the following algorithm.

{\bf Algorithm 3.6.}\\
1) For each fixed $\xi=\overline{r+1,s},$ we solve $IP(e_\xi, G),$
and find the potential $q$ on $e_\xi$.\\
2) For each fixed $k=\overline{1,p-1},$ we solve $IP(e_k, G),$
and find the potential $q$ on $e_k$.\\
3) For $\mu=\om-1,\om-2,\ldots,1,0,$ we successively perform the following
operations: \\
For each fixed simple edge $e_k\in{\cal E}^{(\mu)}$, $p\le k\le r,$
we construct the potential $q$ on $e_k$ by the descent procedure.

\smallskip
{\bf Remark 3.7. } Let $p\le 1$ (i.e. $p=0$ or $p=1$). Then the inverse
problem is formulated as follows: Given the spectra $\Lambda_{\nu}^\xi,\;
\xi=\overline{r+1,s},\; \nu=0,1,$ construct the potential $q$ on $G.$
For $p=1$ all arguments and results remain true; in particular,
for finding the potential $q$ on $G,$ one can use Algorithm 3.6 without Step 2.
If $p=0,\; r>0,$ Then the tree $G^{*}$ is not empty. We choose and fix one
of the boundary vertices of $G^{*}$ as a root, and repeat the above-mentioned
procedure. If $r=0,$ then the tree $G^{*}$ is empty, and we omit Step 3 in
Algorithm 3.6.

\smallskip
{\bf Remark 3.8. } The goal of this remark is to compare the M-matrix from
[53] with the well-known spectral data for the classical Sturm-Liouville
operators on an interval. Let $Z_j(x_j,\la),$ $W_j(x_j,\la)$ be solutions of
equation (3.4) with the initial conditions
$Z_j(l_j,\la)=W'_j(l_j,\la)=1,\; Z'_j(l_j,\la)=W_j(l_j,\la)=0.$ Denote
$$
m_{j,\nu}(\la)=\frac{Z_j^{(\nu)}(l_j/2,\la)}{W_j^{(\nu)}(l_j/2,\la)},\quad
m_{j+s,\nu}(\la)=\frac{C_j^{(\nu)}(l_j/2,\la)}{S_j^{(\nu)}(l_j/2,\la)},\quad
\nu=0,1,\; j=\overline{1,s}.
$$
For each fixed $j=\overline{1,s}$ and $\nu=0,1,$ these functions are the Weyl
functions (see [12]) for equation (3.4) on the intervals $(l_j/2,l_j)$ and
$(0,l_j/2),$ respectively. Fix $j=\overline{1,s}$ and $\nu=0,1.$
It is known (see [12]) that the specification of $m_{j,\nu}(\la)$ and
$m_{j+s,\nu}(\la)$ uniquely determines the potential $q_j$ on $(0,l_j).$
This is the classical inverse Sturm-Liouville problem on an interval.

In [53] equation (3.4) is studied on a graph with some boundary conditions.
The so-called M-matrix is introduced (see [53] for details), and the
uniqueness of recovering equation (3.4) from the M-matrix is considered.
It is easy to check that the specification of the M-matrix from [53]
is equivalent to the specification of the Weyl functions $m_{j,\nu}(\la)$
for all $j=\overline{1,2s}$ and $\nu=0,1.$ This means that in [53] there is
nothing specific for a graph, and the inverse problem is considered really
not on a graph but separately for each fixed $j=\overline{1,s}$ on the finite
interval $(0,l_j).$ In other words, the M-matrix does not reflect spectral properties
of operators on a graph. We also note that the M-matrix is introduced in
[53] as the Weyl-type matrix for a $2s$- dimensional matrix Sturm-Liouville
equation with separated boundary conditions and with a diagonal potential.
The inverse problem for the matrix Sturm-Liouville equation has been solved in
many papers (see [54] and the references therein). Below are two examples
on Remark 3.8.

{\it Examples. } 1) Let $s=2,\,l_1=l_2=1,$ and let the boundary conditions have
the form $y_2(0)=0,\, y_1(0)=y_1(1)=y_2(1),\, y'_1(0)=y'_1(1)+y'_2(1).$
These conditions are related to a graph with one cycle and one boundary edge.
Let $M(\la)=[M_{kr}(\la)]_{k,r=\overline{1,4}}$ be the M-matrix from [53]. Then
$$
M_{11}(\la)-M_{13}(\la)=m_{10}(\la)/2,\quad M_{11}(\la)+M_{13}(\la)=m_{11}(\la)/2,
$$
$$
M_{22}(\la)-M_{24}(\la)=m_{20}(\la)/2,\quad M_{22}(\la)+M_{24}(\la)=m_{21}(\la)/2,
$$
$$
M_{31}(\la)-M_{33}(\la)=m_{30}(\la)/2,\quad M_{31}(\la)+M_{33}(\la)=-m_{31}(\la)/2,
$$
$$
M_{42}(\la)-M_{44}(\la)=-2/m_{40}(\la),\quad M_{42}(\la)+M_{44}(\la)=2/m_{41}(\la).
$$
We do not need other elements of the M-matrix.

2) Let $s=1,\,l_1=1,$ and let the boundary conditions have
the form $y_1(0)-y_1(1)=y'_1(0)-y'_1(1)=0.$
These conditions are related to the classical periodic problem.
Let $M(\la)=[M_{kr}(\la)]_{k,r=\overline{1,2}}$ be the M-matrix from [53]. Then
$$
M_{11}(\la)-M_{12}(\la)=m_{10}(\la)/2,\quad M_{11}(\la)+M_{12}(\la)=m_{11}(\la)/2,
$$
$$
M_{21}(\la)-M_{22}(\la)=-2/m_{20}(\la),\quad M_{21}(\la)+M_{22}(\la)=2/m_{21}(\la).
$$

{\bf Remark 3.9. } An important and wide subclass of graphs is the so-called
A-graphs (another name is cactus-graph), having the property that any two cycles
can have at most one common point. Inverse problems for compact A-graphs were
studied in [17-18], where an algorithm for the solution of the inverse problem
was obtained and the uniqueness of the solution was proved.

\begin{center}
{\bf 4. Arbitrary order differential operators on compact graphs}
\end{center}

In this section we present results on inverse spectral problems for arbitrary order
differential operators on compact graphs. In subsection 4.1 we study $n$-th order
differential operators for a fixed $n>2$ on an arbitrary tree.
In subsection 4.2 we provide the solution of the inverse problem for variable order
differential operators when differential equations have different orders on different
edges. As the main spectral characteristic we introduce and study the so-called
Weyl matrices which are a generalization of the Weyl function for the classical
Sturm-Liouville operator and a generalization of the Weyl matrix for arbitrary order
equations on an interval and on graphs [12]. It is proved that the specification
of the Weyl matrices uniquely determines the coefficients of the differential equation
on the graph. We also provide a procedure for constructing the solution of the
inverse problem. In subsection 4.3 we describe the solution of the inverse problem
for variable order differential operators ona graph with a cycle. The results
of subsection 4.3 were obtained in [50].

\medskip
{\bf 4.1. Inverse problems for $n$-th order operators.}
Consider a compact, connected tree $T$ in
${\bf R^m}$ with the root $v_0$, the set of vertices $V=\{v_0,\ldots,
v_r\}$ and the set of edges ${\cal E}=\{e_1,\ldots, e_r\}$.
We suppose that the length of each edge is equal to $1$.
Without loss of generality we assume that $v_0$ is a boundary vertex.
We will use notations from subsection 2.1 related to the tree $T$.

A function $Y$ on $T$ may be represented as $Y(x)=\{y_j(x)\}_{j\in J}$,
$x\in [0,1],$ where $J:=\{j:\;j=\overline{1,r}\}$, and the function
$y_j(x)$ is defined on the edge $e_j$. Fix $n\ge 2.$ Let $q_\nu(x)=
\{q_{\nu j}(x)\}_{j=\overline{1,r}}$, $\nu=\overline{0,n-2}$ be integrable
complex-valued functions on $T.$ Consider the following $n$-th order
differential equation on $T$:
$$
y^{(n)}_j(x)+\sum_{\nu=0}^{n-2} q_{\nu j}(x)y^{(\nu)}_j(x)
=\la y_j(x),\quad  j=\overline{1,r},                                        \eqno(4.1)
$$
and $y^{(\nu)}_j(x) \in AC[0,1],$ $j=\overline{1,r},$ $\nu=\overline{0,n-1}.$
Denote by $q=\{q_\nu\}_{\nu=\overline{0,n-2}}$ the set of the coefficients
of equation (1); $q$ is called the potential.

For each internal vertex $v_m$, $m=\overline{p+1,r},$ we denote
by $R(v_m):=\{e\in {\cal E}:\;e=[v_m,w],\, w\in V\}$ the set of
edges emanating from $v_m$. Let $R_m:=\{ j:\; e_j\in R(v_m)\},$
and let $\om_m$ be the number of edges emanating from $v_m$.
If $R_m=\{\al_{jm}\}_{j=\overline{1,\om_m}}$, then we put
$y_{(jm)}:=y_{\al_{jm}}$. Consider the linear forms
$$
U_{j\nu m}(Y)=\sum_{\mu=0}^{\nu} \ga_{j\nu\mu m}y_{(jm)}^{(\mu)}(1),
\quad m=\overline{p+1,r},\;j=\overline{1,\om_m},\;\nu=\overline{0,n-1},
$$
where $\ga_{j\nu\mu m}$ are complex numbers, $\ga_{j\nu m}:=
\ga_{j\nu\nu m}\ne 0,$ and satisfy the regularity condition for
matching [44]. Let $\Psi_{sk}(x,\la)=
\{\psi_{skj}(x,\la)\}_{j=\overline{1,r}}$,
$s=\overline{1,p}$, $k=\overline{1,n},$ be solutions of
equation (4.1) satisfying the boundary conditions
$$
\left.\begin{array}{c}
y_s^{(\nu-1)}(0)=\de_{k\nu},\quad \nu=\overline{1,k},\\[3mm]
y_j^{(\xi-1)}(\de_{j,p+1})=0,\quad \xi=\overline{1,n-k},\;
j=\overline{1,p+1}\setminus s,
\end{array}\right\}                                                     \eqno(4.2)
$$
and the matching conditions in each internal vertex
$v_m,\;m=\overline{p+1,r}$:
$$
U_{j\nu m}(Y)+y_m^{(\nu)}(0)=0,\quad
j=\overline{1,\om_m},\;\nu=\overline{0,k-1},                             \eqno(4.3)
$$
$$
\di\sum_{j=1}^{\om_m} U_{j\nu m}(Y)+y_m^{(\nu)}(0)=0,
\quad \nu=\overline{k,n-1}.                                              \eqno(4.4)
$$
Here and in the sequel, $\de_{k\nu}$ is the Kronecker symbol. The solution
$\Psi_{sk}$ is called the Weyl solution of  order $k$ with respect to the
boundary vertex $v_s$. We introduce the matrices
$M_{s}(\la)=[M_{sk\nu}(\la)]_{k,\nu=\overline{1,n}},\; s=\overline{1,p},$
where $M_{sk\nu}(\la):=\psi^{(\nu-1)}_{sks}(0,\la).$ It follows from
the definition of $\psi_{skj}$ that $M_{sk\nu}(\la)=\de_{k\nu}$ for
$k\ge\nu,$  and $\det M_s(\la)\equiv 1.$ The matrix $M_s(\la)$
is called the Weyl matrix with respect to the boundary vertex
$v_s$. Denote by $M=\{M_s\}_{s=\overline{1,p}}$ the set of the
Weyl matrices. The inverse problem is formulated as follows.

{\it Inverse problem 4.1.} Given $M,$ construct $q$ on $T.$

This inverse problem is a generalization of the well-known
inverse problems for differential operators on an interval
and on graphs [12].

Let $\{C_{\mu j}(x,\la)\}_{\mu=\overline{1,n}}$, $j=\overline{1,r}$
be solutions of equation (4.1) on the edge $e_j$ under the
initial conditions $C^{(\nu-1)}_{\mu j}(0,\la)=\de_{\mu\nu}$,
$\mu,\nu=\overline{1,n}.$ One has
$$
\psi_{skj}(x,\la)=
\sum_{\mu=1}^{n} M_{skj\mu}(\la)C_{\mu j}(x,\la),\quad
s=\overline{1,p},\; k=\overline{1,n},\; j=\overline{1,r},        \eqno(4.5)
$$
where the coefficients $M_{skj\mu}(\la)$ do not depend on
$x.$ In particular, $M_{sks\mu}(\la)=M_{sk\mu}(\la),$ and
$$
\psi_{sks}(x,\la)=C_{ks}(x,\la)
+\sum_{\mu=k+1}^{n} M_{sk\mu}(\la)C_{\mu s}(x,\la).              \eqno(4.6)
$$
Substituting the representation (4.5) into (4.2)-(4.4) for
$Y=\Psi_{sk}$, we obtain a linear algebraic system with respect
to $M_{skj\mu}(\la).$ Solving this system by Cramer's rule one gets
$M_{skj\mu}(\la)=\Delta_{skj\mu}(\la)(\Delta_{sk}(\la))^{-1},$
where the functions $\Delta_{skj\mu}(\la)$ and $\Delta_{sk}(\la)$
are entire in $\la$ of order $1/n$. Thus, the functions
$M_{skj\mu}(\la)$ are meromorphic in $\la,$ and consequently,
the Weyl solutions and the Weyl matrices are
meromorphic in $\la.$ In particular,
$$
M_{sk\mu}(\la)=
\Delta_{sk\mu}(\la)(\Delta_{sk}(\la))^{-1},\quad k<\mu,              \eqno(4.7)
$$
where $\Delta_{sk\mu}(\la)=\Delta_{sks\mu}(\la),\;
\Delta_{sk}(\la)=\Delta_{skk}(\la).$ The zeros
$\Lambda_{sk\mu}:=\{\la_{lsk\mu}\}_{l\ge 1}$ of the function
$\Delta_{sk\mu}(\la)$ coincide with the eigenvalues of
the boundary value problem $L_{sk\mu}$ for equation (4.1)
with the matching conditions (4.3)-(4.4) and with the boundary
conditions $y_s^{(\nu-1)}(0)=0,$ $\nu=\overline{1,k-1},\mu,$
$y_j^{(\xi-1)}(\de_{j,p+1})=0,$ $\xi=\overline{1,n-k},\;
j=\overline{1,p+1}\setminus s.$

\smallskip
{\it Problem} $Z_k(T, v_0, \{a_\nu\}_{\nu=\overline{1,k}})$.
Fix $k=\overline{1,n-1}.$ Let $a_\nu$, $\nu=\overline{1,k},$
be complex numbers. Denote by $\vv_k=\{\vv_{kj}\}_{j\in J}$
the solution of equation (4.1) on $T$ satisfying the matching
conditions (4.3)-(4.4) and the following boundary conditions
$$
(\vv^{(\nu-1)}_{k})_{|v_0}=a_{\nu},\;\nu=\overline{1,k},
\quad (\vv^{(\xi-1)}_{k})_{|v_j}=0,\;
v_j\in\Gamma\setminus v_0,\;\xi=\overline{1,n-k}.                   \eqno(4.8)
$$
Denote $m_{kj\mu}(\la)=\vv_{kj}^{(\mu-1)}(0,\la),\;j\in J,
\;\mu=\overline{1,n}.$ Then
$$
\vv_{kj}(x,\la)=\sum_{\mu=1}^{n} m_{kj\mu}(\la)C_{\mu j}(x,\la).   \eqno(4.9)
$$
The matrix $m_{kj\mu}(\la),\; j\in J,\;\mu=\overline{1,n}$ is called
the transition matrix for $T$ with respect to $v_0$ and $\{a_\nu\}.$
Substituting (4.9) into (4.8) and (4.3)-(4.4) for $Y=\vv_k$ we obtain a
linear algebraic system with respect to $m_{kj\mu}(\la),$ $j\in J,\;\mu=
\overline{1,n}.$ The determinant of this system is $\Delta_{sk}(\la).$
Solving this system by Cramer's rule we find the transition matrix
$m_{kj\mu}(\la),\; j\in J,\;\mu=\overline{1,n}$ for $T$ with
respect to $v_0$ and $\{a_\nu\}.$ The algebraic problem of calculating
the transition matrix $m_{kj\mu}(\la),\; j\in J,\;\mu=\overline{1,n}$
by Cramer's rule is called problem $Z_k(T,v_0, \{a_\nu\}_{\nu=
\overline{1,k}}).$

Fix $s=\overline{1,p},$ and consider the following auxiliary
inverse problem on the edge $e_s$, which is called IP(s):
{\it Given $M_s$, construct the functions $q_{\nu s}$,
$\nu=\overline{0,n-2}$ on the edge $e_s$.}

\smallskip
{\bf Theorem 4.2. }{\it The specification of the matrix
$M_s$ uniquely determines the potential on the edge $e_s$.}

\smallskip
This theorem was proved in [44], where a constructive procedure for
the solution of the inverse problem IP(s) was also obtained.
Let us introduce the Weyl matrices for internal vertices. Fix
$v_m\notin\Gamma.$ Denote $T_m^0:=\{z\in T:\; v_m<z\}$,
$T_m:=T\setminus T_m^0$. Clearly, $T_m$ is a rooted tree with the root
$v_0.$ Let $\Gamma_m$ be the set of boundary vertices of $T_m$. Denote
$J_m:=\{j:\; e_j\in T_m\}.$ If $Y=\{y_j\}_{j\in J}$ is a function on
$T$, then $\{Y\}_m:=\{y_j\}_{j\in J_m}$ is a function on $T_m$.

Fix $v_m\notin\Gamma$, $k=\overline{1,n-1}.$ Let
$\Psi_{mk}(x,\la)=\{\psi_{mkj}(x,\la)\}_{j\in J_m}$ be
the solution of equation (4.1) on $T_m$ satisfying (4.3)-(4.4) on
$T_m$ and the boundary conditions
$(\Psi^{(\nu-1)}_{mk})_{|v_m}=\de_{\nu k},$ $\nu=\overline{1,k},$
$(\Psi^{(\xi-1)}_{mk})_{|v_j}=0,$ $v_j\in\Gamma_m\setminus v_m,
\;\xi=\overline{1,n-k}.$ The function $\Psi_{mk}$ is the Weyl
solution of (4.1) on $T_m$ with respect to the internal vertex
$v_m$. We introduce the matrix
$M_{m}(\la)=[M_{mk\nu}(\la)]_{k,\nu=\overline{1,n}},$
where $M_{mk\nu}(\la):=\psi^{(\nu-1)}_{mkm}(0,\la).$ Then
$M_{mk\nu}(\la)=\de_{k\nu}$ for $k\ge\nu,$  and
$\det M_m(\la)\equiv 1.$ The matrix $M_m(\la)$
is called the Weyl matrix for $T_m$ with respect to
the internal vertex $v_m$.

\smallskip
{\bf Theorem 4.3. }{\it Fix $v_m\notin\Gamma$ and
$k=\overline{1,n-1}.$ Let $e_s=[v_m,v_s]\in R(v_m).$ Then}
$$
M_{mk\nu}(\la)=\di\frac{\det[\psi_{s\mu m}(0,\la),\ldots,
\psi_{s\mu m}^{(k-2)}(0,\la),
\psi_{s\mu m}^{(\nu-1)}(0,\la)]_{\mu=\overline{1,k}}}
{\det[\psi_{s\mu m}^{(\xi-1)}(0,\la)]_{\xi,\mu=
\overline{1,k}}}\,,\;\;\nu=\overline{k+1,n}.                       \eqno(4.10)
$$
Using the fundamental system of solutions
$\{C_{\mu j}(x,\la)\}_{\mu=\overline{1,n}}$, one gets that
(4.5) holds for $s=\overline{1,r},\; k=\overline{1,n},\; j\in J_s$,
where $J_s:=J$ for $s=\overline{1,p},$ and the coefficients
$M_{skj\mu}(\la)$ do not depend on $x.$ In particular,
$M_{sks\mu}(\la)=M_{sk\mu}(\la),$ and (4.6) holds for
$s=\overline{1,r},\; k=\overline{1,n}$. Therefore, we obtain
$$
\psi_{skj}^{(\nu-1)}(1,\la)=\sum_{\mu=1}^{n}
M_{skj\mu}(\la)C_{\mu j}^{(\nu-1)}(1,\la),\quad
k,\nu=\overline{1,n},\; s=\overline{1,r},\; j\in J_s,             \eqno(4.11)
$$
$$
\psi_{sks}^{(\nu-1)}(1,\la)=C_{ks}^{(\nu-1)}(1,\la)+
\sum_{\mu=k+1}^{n} M_{sk\mu}(\la)C_{\mu s}^{(\nu-1)}(1,\la),
\quad k,\nu=\overline{1,n},\; s=\overline{1,r}.                   \eqno(4.12)
$$

Now we will obtain a constructive procedure for
the solution of inverse problem 4.1 and prove its uniqueness. Let
the set $M=\{M_s\}_{s=\overline{1,p}}$ for equation (4.1) be given.
The procedure for the solution of inverse problem 1 consists in
the realization of the so-called $A_\xi$- procedures successively
for $\xi=\sigma,\sigma-1,\ldots,1,$ where $\sigma$ is the height
of the tree $T.$ Let us describe $A_\xi$- procedures by induction.
Fix $\xi=\overline{1,\sigma},$ and suppose that
$A_\sigma,\ldots,A_{\xi+1}$- procedures have been already
carried out. Let us carry out $A_\xi$- procedure.

{\bf ${\bf A_\xi}$- procedure. } For each $v_s\in V^{(\xi)},$
the Weyl matrix $M_s(\la)$ is given. Indeed, if
$v_s\in V^{(\xi)}\cap\Gamma,$ then $M_s(\la)$ is given a priori,
and if $v_s\in V^{(\xi)}\setminus\Gamma,$ then $M_s(\la)$ was
calculated on the previous steps according to
$A_\sigma,\ldots,A_{\xi+1}$- procedures.

1) For each edge $e_s\in{\cal E}^{(\xi)}$ we solve the local inverse
problem IP(s) and find the coefficients $q_{\nu s}(x),$ $x\in[0,1],
\; \nu=\overline{0,n-2}$ of equation (4.1) on the edge $e_s$.
If $\xi=1,$ then inverse problem 4.1 is solved, and we stop our
calculations. If $\xi>1,$ we go on to the next step.

2) For each edge $e_s\in{\cal E}^{(\xi)},$ we construct
$C_{\mu s}(x,\la),$ $x\in[0,1],\;\mu=\overline{1,n},$ and
calculate the functions
$$
\psi^{(\nu-1)}_{sks}(1,\la),\quad k,\nu=\overline{1,n}.           \eqno(4.13)
$$

3) {\it Returning procedure.} For each fixed $v_m\in V^{(\xi-1)}
\setminus\Gamma,$ $k=\overline{1,n},$ and for all $e_i,e_s\in R(v_m),
\,i\ne s,$ we fulfill the following operations:

(i) Using (4.13) and the matching conditions (4.3) in $v_m$ for
$Y=\Psi_{sk}$, we obtain
$$
\psi^{(\nu-1)}_{ski}(1,\la)=a_{ski\nu}(\la),\quad \nu=\overline{1,k},
$$
where $a_{ski\nu}(\la)$ are constructed as linear
combinations of the functions (4.13) for $\nu=\overline{1,k}$.

(ii) Consider the tree $T_i^1:=T_i^0\cup\{e_i\}$ with the
root $v_m$. Solving the algebraic problem
$Z_k(T^1_i,v_m,\{a_{ski\nu}(\la)\}_{\nu=\overline{1,k}}),$
we calculate the transition matrix $\{M_{skj\mu}(\la)\},$
$\mu=\overline{1,n},\; j\in J_{i1}$, where $J_{i1}:=
\{j:\;e_j\in T_i^1\}.$

4) For each fixed $v_m\in V^{(\xi-1)}\setminus\Gamma,$ and
for all $e_j,e_s\in R(v_m),\,j\ne s,$ we construct the functions
$$
\psi^{(\nu-1)}_{skj}(1,\la),\quad k,\nu=\overline{1,n},          \eqno(4.14)
$$
by (4.11). Furthermore, using (4.13), (4.14) and the matching
conditions (4.3)-(4.4) we find the functions
$\psi^{(\nu-1)}_{skm}(0,\la),$ $k,\nu=\overline{1,n}.$

5) For each fixed $v_m\in V^{(\xi-1)}\setminus\Gamma,$
we calculate the matrix $M_m(\la)$ via (4.10).

Thus, we have obtained the solution of inverse problem 4.1
and proved its uniqueness, i.e. the following assertion holds.

\smallskip
{\bf Theorem 4.4 }{\it The specification of the Weyl matrices
$M=\{M_s\}_{s=\overline{1,p}}$ uniquely determines the potential
$q$ on $T.$ The solution of inverse problem 4.1 can be constructed
by executing successively $A_\sigma,A_{\sigma-1},\ldots,A_{1}$-
procedures.}

\smallskip
Let us now obtain the solution of the inverse problem
of recovering equation (4.1) from discrete spectral characteristics.
Denote by $\Lambda:=\{\Lambda_{sk\mu}\},$ $s=\overline{1,p},$
$k=\overline{1,n-1},$ $\mu=\overline{k,n},$ the set of spectra
of the boundary value problems $L_{sk\mu}$. The inverse problem
is formulated as follows.

{\it Inverse problem 4.5 } Given the set of spectra $\Lambda,$
construct $q$ on $T.$

This inverse problem is a generalization of the classical inverse problem
of recovering Sturm-Liouville operators on an interval from two spectra
and a generalization of the inverse problem of recovering higher-order
differential operators on an interval from spectra (see [12]). The
solution of inverse problem 4.5 can be reduced to the solution of
Inverse problem 4.1. Indeed, using (4.7) and Hadamard's theorem, one gets
$$
M_{sk\mu}(\la)=\be_{sk\mu}\prod_{l=1}^{\iy}
\Big(1-\frac{\la}{\la_{lsk\mu}}\Big)
\Big(1-\frac{\la}{\la_{lsk}}\Big)^{-1},\quad
s=\overline{1,p},\;k=\overline{1,n-1},\; \mu=\overline{k+1,n},         \eqno(4.15)
$$
where $\la_{lsk}:=\la_{lskk}$, and the constants $\be_{sk\mu}$
can be uniquely calculated from the asymptotics of
$M_{sk\mu}(\la)$.
Thus, using the given spectra $\Lambda,$ we can construct uniquely
the set $M=\{M_s\}_{s=\overline{1,p}}$ of the Weyl matrices.
The following assertion holds.

\smallskip
{\bf Theorem 4.6. }{\it The specification of the set of spectra
$\Lambda$ uniquely determines the potential $q$ on $T.$ The solution
of inverse problem 4.5 can be obtained by two steps:

1) Construct the Weyl matrices $M=\{M_s\}_{s=\overline{1,p}}$
by (4.15);

2) Find $q$ on $T$ by executing successively
$A_\sigma,A_{\sigma-1},\ldots,A_{1}$- procedures.}

\smallskip
{\bf Remark 4.7. }{\it Inverse problem from the spectral data.}
It is also possible to consider an inverse problem from other
discrete spectral characteristics, namely from the so-called
spectral data, which describe main parts of Laurent's series
for the Weyl matrices in the neighborhoods of their poles.
For simplicity, let all poles be simple. Denote $\al_{lsk\mu}=
\mathop{\rm Res}\limits_{\la=\la_{lsk}} M_{sk\mu}(\la).$
The data $\Lambda_1:=\{\la_{lsk},\al_{lsk\mu}\},$ $l\ge 1,\,
s=\overline{1,p},\,k=\overline{1,n-1},\,\mu=\overline{k+1,n},$
are called the spectral data. The inverse problem is formulated
as follows: Given the spectral data $\Lambda_1$, construct
the potential $q$ on $T.$
This inverse problem is a generalization of the classical
inverse problem of recovering Sturm-Liouville operators on an
interval (cv/ [9-10]). One can get
that the specification of $\Lambda_1$ uniquely determines $M.$
Thus, the solution of the inverse problem from the spectral data
can be reduced to the solution of inverse problem 4.1.

\medskip
{\bf 4.2. Inverse problems for variable order operators on graphs without cycles.}
We study inverse spectral problems for ordinary differential operators of
variable order on a compact graph without cycles, when differential equations
have different orders on different edges. For simplicity we confine ourselves
to the case of star-graphs. The general case is treated analogously.

Consider a compact star-type graph $T$ in ${\bf R^\om}$ with the set of vertices
$V=\{v_0,\ldots, v_p\}$ and the set of edges ${\cal E}=\{e_1,\ldots, e_p\},$ where
$v_1,\ldots, v_{p}$ are the boundary vertices, $v_0$ is the internal vertex, and
$e_j=[v_{j},v_0],$ $e_1\cap\ldots\cap e_p=\{v_0\}$.
Let $l_j$ be the length of the edge $e_j$. Each edge $e_j\in {\cal E}$ is
parameterized by the parameter $x_j\in [0,l_j]$ such that $x_j=0$ corresponds
to the boundary vertices $v_1,\ldots, v_{p}$, and $x_j=l_j$ corresponds to the
internal vertex $v_0$. A function $Y$ on $T$ may be represented as
$Y=\{y_j\}_{j=\overline{1,p}}$, where the function $y_j(x_j)$ is defined on
the edge $e_j$.
Fix $m=\overline{1,p}.$ Let $n_i, p_i$, $i=\overline{1,m},$ be positive integers
such that $n_1>n_2>\ldots>n_m>1,$ $0<p_1<p_2<\ldots<p_{m-1}<p_m:=p,$ и положим
$n_{m+1}:=1,\; p_0:=0.$ Consider the differential equations on $T$:
$$
y_j^{(n_i)}(x_j)+\sum_{\mu=0}^{n_i-2}q_{\mu j}(x_j)y_j^{(\mu)}(x_j)
=\la y_j(x_j),\quad x_j\in (0, l_j),
\; i=\overline{1,m},\;j=\overline{p_{i-1}+1,p_i},                            \eqno(4.16)
$$
where $q_{\mu j}(x_j)$ are complex-valued integrable functions. Thus,
the differential equations have order $n_i$ on the edges
$e_j,\;j=\overline{p_{i-1}+1,p_i}$. We call $q_j=\{q_{\mu j}\}$ the
potential on the edge $e_j$, and we call $q=\{q_{j}\}_{j=\overline{1,p}}$
the potential on the graph $T.$
Denote $\mu_i=p_i-p_{i-1},\; i=\overline{1,m}.$

Fix $i=\overline{1,m},\; j=\overline{p_{i-1}+1,p_i}.$ Let $\{C_{kj}(x_j,\la)\},$
$k=\overline{1,n_i}$, be the fundamental system of solutions of equation (4.16) on
the edge $e_j$ under the initial conditions $C^{(\mu-1)}_{kj}(0,\la)=\de_{k\mu}$,
$k,\mu=\overline{1,n_i}$. Here and in the sequel, $\de_{k\mu}$ is the Kronecker
symbol. Consider the linear forms
$$
U_{j\nu}(y_j)=\sum_{\mu=0}^{\nu}\ga_{j\nu\mu}y_j^{(\mu)}(l_j),
\quad j=\overline{1,p},
$$
where $\ga_{j\nu\mu}$ are complex numbers, $\ga_{j\nu}:=\ga_{j\nu\nu}\ne 0,$
$\nu=\overline{0,n_i-1}$ for $j=\overline{p_{i-1}+1,p_i}$.

Denote $\langle n\rangle:=(|n|+n)/2,$ i.e. $\langle n\rangle=n$
for $n\ge 0$, and $\langle n\rangle=0$ for $n\le 0.$ Fix
$i=\overline{1,m},\;s=\overline{p_{i-1}+1,p_i},\;k=\overline{1,n_i-1}.$
We introduce the solutions $\Psi_{sk}=\{\psi_{skj}\}_{j=\overline{1,p}}$
of Eq. (4.16) on the graph $T$ as follows. Let $\xi=\overline{i,m};\;
k=\overline{n_{\xi+1},n_\xi-1}.$ Then $\Psi_{sk}$ satisfies the boundary
conditions
$$
\psi_{sks}^{(\eta-1)}(0)=\de_{k\eta},\quad \eta=\overline{1,k},          \eqno(4.17)
$$
$$
\psi_{skj}^{(r)}(0)=0,\quad r=\overline{0,\langle n_l-k-1\rangle};
\; l=\overline{1,m},\; j=\overline{p_{l-1}+1,p_l},\; j\ne s,             \eqno(4.18)
$$
and the matching conditions at the vertex $v_0$:
$$
U_{p_l,\nu}(\psi_{sk,p_l})=U_{j\nu}(\psi_{skj}),\; l=\overline{\xi,m},
\;j=\overline{1,p_l-1},\; \nu=\overline{n_{l+1}-1,\min(k-1,n_l-2)},      \eqno(4.19)
$$
$$
\di\sum_{j=1}^{p_\xi} U_{j\nu}(\psi_{skj})=0,\; \nu=\overline{k,n_\xi-1};
\quad \di\sum_{j=1}^{p_l} U_{j\nu}(\psi_{skj})=0,\; l=\xi-1,\ldots, i,
\; \nu=\overline{n_{l+1},n_l-1}.                                         \eqno(4.20)
$$
The function $\Psi_{sk}$ is called the Weyl solution of order $k$ with
respect to the boundary vertex $v_s$. Additionally we define
$\psi_{sn_is}(x_s,\la):=C_{n_is}(x_s,\la).$

Let us introduce the matrices $M_{s}(\la),\;s=\overline{p_{i-1}+1,p_i},
\;i=\overline{1,m}$:
$$
M_{s}(\la)=[M_{sk\mu}(\la)]_{k,\mu=\overline{1,n_i}},\quad
M_{sk\mu}(\la):=\psi^{(\mu-1)}_{sks}(0,\la).
$$
It follows from the definition of $\Psi_{sk}$ that $M_{sk\mu}(\la)=
\de_{k\mu}$ for $k\ge\mu.$ The matrix $M_s(\la)$ is called the Weyl matrix
with respect to the vertex $v_s$. Fix $N=\overline{1,m}.$

\smallskip
{\bf Inverse problem 4.8.} Given $\{M_{s}(\la)\},\;
s=\overline{1,p}\setminus p_N$, construct $q$ on $T.$

\smallskip
Fix $i=\overline{1,m},\; s=\overline{p_{i-1}+1,p_i}.$
It follows from the boundary conditions (4.17) for the Weyl solutions that
$$
\psi_{sks}(x_s,\la)=C_{ks}(x_s,\la)+\di\sum_{\mu=k+1}^{n_i}
M_{sk\mu}(\la)C_{\mu s}(x_s,\la), \quad k=\overline{1,n_i}.                \eqno(4.21)
$$
Using the fundamental system of solutions $\{C_{\mu j}(x_j,\la)\}$
on the edge $e_j$, we get
$$
\psi_{skj}(x_j,\la)=\sum_{\mu=1}^{n_l}M_{skj\mu}(\la)C_{\mu j}(x_j,\la),
\;j=\overline{p_{l-1}+1,p_l},\;l=\overline{1,m},\;k=\overline{1,n_i-1},     \eqno(4.22)
$$
where the coefficients $M_{skj\mu}(\la)$ do not depend on $x_j.$ In particular,
$M_{sks\mu}(\la)=M_{sk\mu}(\la).$ Substituting (4.22) into boundary and
matching conditions (4.17)-(4.20) for the Weyl-type solu\-tions $\Psi_{sk},$
we obtain a linear algebraic system $D_{sk}^1(\la)$ with respect to
$M_{skj\mu}(\la).$ Solving this system by Cramer's rule one gets
$M_{skj\mu}(\la)=\Delta_{skj\mu}(\la)/\Delta_{sk}(\la),$
where the functions $\Delta_{skj\mu}(\la)$ and $\Delta_{sk}(\la)$
are entire in $\la.$ In particular,
$M_{sk\mu}(\la)=\Delta_{sk\mu}(\la)/\Delta_{sk}(\la),\; k<\mu,$
where $\Delta_{sk\mu}(\la):=\Delta_{sks\mu}(\la).$ The functions
$M_{skj\mu}(\la)$ are meromorphic in $\la.$

Fix $s=\overline{1,p},$ and consider the following inverse problem
on the edge $e_s$.

\smallskip
{\bf Inverse problem 4.9. } Given the Weyl matrix $M_s$,
construct the potential $q_{s}$ on the edge $e_s$.

\smallskip
This inverse problem was solved in [48], where the uniqueness theorem was
proved and an algorithm for the solution of Inverse problem 4.9 was obtained.

\smallskip
Fix $i=\overline{1,m},\;j=\overline{p_{i-1}+1,p_i}.$ Let
$\vv_{jk}(x_j,\la),\;k=\overline{1,n_i}$ be solutions of Eq. (4.16) on the
edge $e_j$ with the conditions $\vv_{jk}^{(\nu-1)}(l_j,\la)=\de_{k\nu},
\;\nu=\overline{1,k},\;\vv_{jk}^{(\mu-1)}(0,\la)=0,\;
\mu=\overline{1,n_i-k}.$ We introduce the matrix
$m_j(\la)=[m_{jk\nu}(\la)]_{k,\nu=\overline{1,n_i}},$
where $m_{jk\nu}(\la):=\vv^{(\nu-1)}_{jk}(l_j,\la).$ The matrix $m_j(\la)$
is a classical Weyl matrix on the edge $e_j$.

\smallskip
{\bf Inverse problem 4.10. } Fix $j=\overline{1,p}.$ Given the matrix
$m_j$, construct the potential $q_{j}$ on the edge $e_j$.

This inverse problem is the classical one, since it is the inverse problem
of recovering a higher-order differential equation on a finite interval
from its Weyl matrix. This inverse problem has been solved in [11],
where the uniqueness theorem for this inverse problem is proved. Moreover,
in [11] an algorithm for the solution of Inverse problem 4.10 is given, and
necessary and sufficient conditions for the solvability of this inverse
problem are provided.

\smallskip
Fix $i=\overline{1,m},\;j=\overline{p_{i-1}+1,p_i}.$
Then for each fixed $s=\overline{1,p_1}\setminus j,$
$$
m_{j1\nu}(\la)=
\frac{\psi_{s1j}^{(\nu-1)}(l_j,\la)}{\psi_{s1j}(l_j,\la)},
\quad \nu=\overline{2,n_i},                                          \eqno(4.23)
$$
$$
m_{jk\nu}(\la)=\di\frac{\det[\psi_{s\mu j}(l_j,\la),\ldots,
\psi_{s\mu j}^{(k-2)}(l_j,\la), \psi_{s\mu j}^{(\nu-1)}
(l_j,\la)]_{\mu=\overline{1,k}}}{\det[\psi_{s\mu j}^{(\xi-1)}
(l_j,\la)]_{\xi,\mu=\overline{1,k}}}\,,\;2\le k<\nu\le n_i.          \eqno(4.24)
$$

\smallskip
{\it Solution of Inverse problem 4.8. } Let us construct an algorithm
for the solution of Inverse problem 4.8.

{\it Step 1. } Let the Weyl-type matrices $\{M_{s}(\la)\},\;
s=\overline{1,p}\setminus p_N$, be given. Solving Inverse problem 4.9
for each fixed $s=\overline{1,p}\setminus p_N,$ we find the potentials
$q_{s}$ on the edges $e_s$, $s=\overline{1,p}\setminus p_N$.

{\it Step 2. } Using the knowledge of the potential on the edges
$e_s$, $s=\overline{1,p}\setminus p_N$, we construct the Weyl-type
matrix $m_{p_N}$.

{\it Step 3. } Solving Inverse problem 4.10 for $j=p_N$ we find
the potential $q_{p_N}$ on $e_{p_N}$.

\smallskip
Steps 1 and 3 have been already studied. It remains to fulfil Step 2.

Suppose that Step 1 was already made, and we found the potentials
$q_{s}$, $s=\overline{1,p}\setminus p_N$, on the edges $e_s$,
$s=\overline{1,p}\setminus p_N$. Then we calculate the functions
$C_{kj}(x_j,\la),$ $j=\overline{1,p}\setminus p_N;$
here $k=\overline{1,n_i}$ for $j=\overline{p_{i-1}+1,p_i}.$

Fix $s=\overline{1,p_1}$ (if $N>1$), and $s=\overline{1,p_1-1}$ (if $N=1$).
All calculations below will be made for this fixed $s.$

Our goal now is to construct the Weyl-type matrix $m_{p_N}(\la).$
By virtue of (4.23)-(4.24), we have to calculate the functions
$$
\psi_{skp_N}^{(\nu)}(l_{p_N},\la),
\quad k=\overline{1,n_N-1},\;\nu=\overline{0,n_N-1}.             \eqno(4.25)
$$
We will find the functions (4.25) by the following steps.

\smallskip
1) Using (4.21) we construct the functions
$$
\psi_{sks}^{(\nu)}(l_s,\la),
\;k=\overline{1,n_N-1},\;\nu=\overline{0,n_1-1},                    \eqno(4.26)
$$
by the formula
$$
\psi_{sks}^{(\nu)}(l_s,\la)=C_{ks}^{(\nu)}(l_s,\la)+
\sum_{\mu=k+1}^{n_1}M_{sk\mu}(\la)C_{\mu s}^{(\nu)}(l_s,\la).        \eqno(4.27)
$$

2) Consider a part of the matching conditions (4.19) on $\Psi_{sk}$. More
precisely, let $\xi=\overline{N,m},\;k=\overline{n_{\xi+1},n_\xi-1},\;
l=\overline{\xi,m},\;j=\overline{1,p_l-1}.$ Then, in particular,
(4.19) yields
$$
U_{p_l,\nu}(\psi_{skp_l})=U_{j\nu}(\psi_{skj}), \quad
\nu=\overline{n_{l+1}-1,\min(k-1,n_l-2)}.                            \eqno(4.28)
$$
Since the functions (4.26) are known, it follows from (4.28) that
one can calculate the functions
$$
\psi_{skj}^{(\nu)}(l_j,\la),\;\xi=\overline{N,m},\;
k=\overline{n_{\xi+1},n_\xi-1},\; l=\overline{\xi,m},\;
j=\overline{1,p_l},\;\nu=\overline{n_{l+1}-1,\min(k-1,n_l-2)}.      \eqno(4.29)
$$
In particular we found the functions (4.25) for $\nu=\overline{0,k-1}.$

\smallskip
3) It follows from (4.22) and the boundary conditions on $\Psi_{sk}$ that
$$
\psi_{skj}^{(\nu)}(l_j,\la)=\di\sum_{\mu=\max(n_l-k+1,2)}^{n_l}
M_{skj\mu}(\la)C_{\mu j}^{(\nu)}(l_j,\la),                       \eqno(4.30)
$$
$$
k=\overline{1,n_1-1},\;l=\overline{1,m},\;
j=\overline{p_{l-1}+1,p_l}\setminus s,\;\nu=\overline{0,n_l-1}.
$$
We consider only a part of relations (4.30). More precisely, let
$\xi=\overline{N,m},\;k=\overline{n_{\xi+1},n_\xi-1},$
$l=\overline{1,m},\;j=\overline{p_{l-1}+1,p_l},\; j\ne p_N,\;
j\ne s,\; \nu=\overline{0,\min(k-1,n_l-2)}.$ Then
$$
\di\sum_{\mu=\max(n_l-k+1,2)}^{n_l}M_{skj\mu}(\la)
C_{\mu j}^{(\nu)}(l_j,\la)=\psi_{skj}^{(\nu)}(l_j,\la),
\quad \nu=\overline{0,\min(k-1,n_l-2)}.                         \eqno(4.31)
$$
For this choice of parameters, the right-hand sides in (4.31)
are known, since the functions (4.29) are known. Relations (4.31)
form a linear algebraic system $\sigma_{skj}$ with respect to the
coefficients $M_{skj\mu}(\la).$ Solving the system by Cramer's
rule we find the functions $M_{skj\mu}(\la).$
Substituting them into (4.30), we calculate the functions
$$
\psi_{skj}^{(\nu)}(l_j,\la),\quad k=\overline{1,n_N-1},\;
l=\overline{1,m},\;j=\overline{p_{l-1}+1,p_l}\setminus p_N,
\; \nu=\overline{0,n_l-1}.                                     \eqno(4.32)
$$
Note that for $j=s$ these functions were found earlier.

\smallskip
4) Let us now use the generalized Kirchhoff's conditions (4.20)
for $\Psi_{sk}$. Since the functions (4.32) are known, one can
construct by (4.20) the functions (4.25) for $k=\overline{1,n_N-1},
\;\nu=\overline{k,n_N-1}.$ Thus, the functions (4.25) are known
for $k=\overline{1,n_N-1},\; \nu=\overline{0,n_N-1}.$

\smallskip
Since the functions (4.25) are known, we construct the Weyl-type
matrix $m_{p_N}(\la)$ via (4.23)-(4.24) for $j=p_N.$ Thus,
we have obtained the solution of Inverse problem 4.8 and proved
its uniqueness, i.e. the following assertion holds.

\smallskip
{\bf Theorem 4.11. }{\it The specification of the Weyl-type matrices
$M_s(\la),$ $s=\overline{1,p}\setminus p_N$, uniquely determines the
potential $q$ on $T.$ The solution of Inverse problem 4.8 can be
obtained by the following algorithm.}

{\bf Алгоритм 4.12. }{\it Given the Weyl matrices $M_s(\la),$
$s=\overline{1,p}\setminus p_N$.

1) Find  $q_{s}$, $s=\overline{1,p}\setminus p_N$, by solving
Inverse problem 4.9 for each fixed $s=\overline{1,p}\setminus p_N$.

2) Calculate $C_{kj}^{(\nu)}(l_j,\la),\; j=\overline{1,p}\setminus p_N$; here
$k=\overline{1,n_i},\;\nu=\overline{0,n_i-1}$ for $j=\overline{p_{i-1}+1,p_i}.$

3) Fix $s=\overline{1,p_1}$ (if $N>1$), and $s=\overline{1,p_1-1}$ (if $N=1$).
All calculations below will be made for this fixed $s.$
Construct the functions (4.26) via (4.27).

4) Calculate the functions (4.29) using (4.28).

5) Find the functions $M_{skj\mu}(\la),$ by solving the linear algebraic
systems $\sigma_{skj}$.

6) Construct the functions (4.25) using (4.20).

7) Calculate the Weyl matrix $m_{p_N}(\la)$ via (4.23)-(4.24) for $j=p_N$.

8) Construct the potential $q_{p_N}$ on the edge $e_{p_N}$
by solving Inverse problem 4.10.}

\medskip
{\bf 4.3. Inverse problems for variable order operators on graphs with cycles.}
Consider a compact graph $G$ with vertices $V = \{ v_0, \dots, v_m \}$
and edges ${\cal E}=\{ e_0, \dots, e_m\}$, where $e_j = [v_j, v_0]$,
$j=\overline{1, m}$, and $e_0$ is a cycle having only one vertex $v_0$.
Thus, $v_j$, $j = \overline{1, m}$ are boundary vertices, and $v_0$ is
a interior vertex. Let $T_j$ be the length of the edge $e_j$. For each
edge $e_j \in{\cal E}$ we introduce the parameter $x_j \in [0, T_j]$
such that for $j = \overline{1, m}$ the endpoint $x_j = 0$ corresponds
to the vertex $v_j$, and the endpoint $x_j=T_j$ corresponds to the vertex
$v_0$. For $j = 0$ both endpoints coincide with $v_0$.

Fix numbers $2=n_0 \le n_1 \le \dots \le n_m$, and consider
differential equations
$$
y_j^{(n_j)}+\sum_{\mu=0}^{n_j-2} q_{\mu j}(x_j) y_j^{(\mu)}(x_j)
=\lambda y_j(x_j), \quad j = \overline{0, m},                            \eqno(4.33)
$$
where $q_{\mu j}\in L[0, T_j]$. The set $q:=\{q_{\mu j}\}_{j=\overline{0,m},
\mu=\overline{0, n_j-2}}$ is called the potential on the graph $G.$
Let us introduce matching conditions at the interior vertex $v_0$, which
are a generalization of Kirchhoff's conditions for Sturm-Liouville
operators and higher order operators on graphs. For this purpose we
consider linear forms
$$
U_{j\nu}(y_j)=\sum_{\mu=0}^{\nu} \gamma_{j\nu\mu}y_j^{(\mu)}(T_j),
\; \gamma_{j\nu\nu}\ne 0,\; j=\overline{1,m},\; \nu=\overline{0, n_j-1},
\qquad U_{0\nu}(y_0)=y_0^{(\nu)}(T_0),\; \nu=0,1,
$$
where $\gamma_{j\nu\mu}$ are complex numbers. We define continuity
conditions $C(\nu)$, $C(0, \alpha)$ and Kirchhoff's conditions
$\mbox{K}(\nu)$ of order $\nu$ as follows:
$$
C(\nu):\quad U_{m\nu}(y_m)=U_{j\nu}(y_j),\; j=\overline{0,m-1},\; \nu<n_j-1,
$$
$$
C(0,\alpha):\quad C(0)\; \mbox{and}\; \alpha y_0(0)=y_0(T_0),\qquad
K(\nu):\quad \sum_{j \colon\nu<n_j} y_j^{(\nu)}(T_j)=\delta_{1\nu} y_0'(0),
$$
where $\delta_{jk}$ is the Kronecker symbol, and $\alpha\ne 0$ is
a complex number.

Let $m>1$ (the case $m=1$ requires small modifications). Fix
$s=\overline{1, m},$ $k=\overline{1, n_s-1}$ and $\mu=\overline{k, n_s}$.
Consider the boundary value problems $L_{sk\mu}$ for system (4.33) with
the boundary conditions
$$
y_k^{(\nu-1)}(0)=0,\quad \nu = \overline{1, k-1}, \mu,
$$
$$
y_j^{(\nu-1)}(0)=0,\quad \nu=\overline{1, n_j-k},\; j=
\overline{1, m}\backslash s:\; n_j > k,
$$
$$
y_j(0)=0,\quad j=\overline{1,m}:\; n_j \le k,
$$
and the matching conditions $C(0, \alpha_s)$, $C(\nu)$, $\nu=\overline{1,k-1}$,
$K(\nu)$, $\nu=\overline{k, n_s-1}$ at the vertex $v_0$. Here $\alpha_s$,
$s=\overline{1, m}$ are nonzero numbers such that among them there are at least
two different numbers. Assume that regularity conditions for matching is
fulfilled. Problems $L_{sk\mu}$ have discrete spectra $\Lambda_{sk\mu} =
\{ \lambda_{lsk\mu}\}_{l\ge 1}$. The characteristic functions $\Delta_{sk\mu}$
of the boundary value problems $L_{sk\mu}$ can be uniquely constructed from
their spectra.

\smallskip
{\bf Inverse problem 4.13. } Given the spectra $\Lambda_{s k \mu}$,
$s=\overline{1,m}$, $k=\overline{1,n_s-1}$, $\mu=\overline{k,n_s}$,
construct the potential $q$ on the graph $G$.

\smallskip
Fix $s=\overline{1,m}$ and $k=\overline{1, n_s-1}$. Let
$\Psi_{sk}=\{\psi_{skj}\}_{j=\overline{1,m}}$ be the solution of system
(4.33) under the conditions
$$
\psi_{sks}^{(\nu-1)}(0)=\delta_{k\nu},\; \nu=\overline{1,k},
$$
$$
\psi_{skj}^{(\xi-1)}(0)=0,\; \xi=\overline{1,n_j-k},\;
j=\overline{1,m}\backslash s \colon k < n_j,
$$
$$
\psi_{skj}(0)=0,\; j=\overline{1,m}\colon k\ge n_j,
$$
$$
C(0, \alpha_s),\; C(\nu),\; \nu=\overline{1,k-1},\;
K(\nu),\;\nu=\overline{k,n_s-1}.
$$
Let $M_{sk\mu}(\lambda):=\psi_{sks}^{(\mu-1)}(0, \lambda)$,
$M_{s n_s \mu}(\lambda)=\delta_{n_s, \mu}$.
The matrix $M_s(\lambda):=[M_{sk\mu}(\lambda)]_{k,\mu=1}^{n_s}$ is called
{\it the Weyl matrix} for the vertex $v_s$. The following relations
hold $M_{sk\mu}(\lambda)=-\Delta_{sk\mu}(\lambda)/\Delta_{skk}(\lambda)$,
i.e. the Weyl matrices can be constructed from the given spectra.
Consider the following auxiliary inverse problem IP(s).

\smallskip
{\bf IP(s).} Given the matrix $M_s$, construct the potential
$q_s:=\{ q_{\mu s}\}_{\mu=0}^{n_s-2}$ on the edge $e_s$.

\smallskip
Problem IP(s) has a unique solution which can be constructed by the method
of spectral mappings (see [11]).

Let $S_0(x_0, \lambda)$ and $C_0(x_0, \lambda)$ be solutions of Eq. (4.33)
on the edge $e_0$ ($n_0=2$) satisfying the initial conditions
$S_0(0, \lambda)=C'_0(0, \lambda)=0$, $S'_0(0, \lambda)=C_0(0, \lambda)=1.$
Using the known characteristic functions $\Delta_{sk\mu}(\lambda)$ and
the potentials $q_s$ on the boundary edges, one can construct
$S_0(T_0,\lambda)$ and $d_{\alpha_s}(\lambda):=C_0(T_0,\lambda)+
\alpha_s S'_0(T_0, \lambda)-\alpha_s-1$. From the functions
$d_{\alpha_s}(\lambda)$ with different $\alpha_s$ it is easy to find
$S'_0(T_0, \lambda)$. It remains to solve on the edge $e_0$ the classical
inverse problem from two spectra (zeros of $S_0(T_0, \lambda)$ and
$S'_0(T_0, \lambda)$), and find the potential $q_{00}(x)$ on the cycle.
Thus, we have obtained the following result.

\smallskip
{\bf Theorem 4.14.} {\it Inverse problem 4.13 is uniquely solvable, and
its solution can be found by the above-mentioned method.}

\begin{center}
{\bf 5. Differential operators on noncompact graphs.}
\end{center}

{\bf 5.1. Inverse Sturm-Liouville problems on a noncompact star-graph.}
Sturm-Liouville differential operators are considered on noncompact star-type
graphs. The graph contains both compact and noncompact edges; this makes
the investigation of inverse problems essentially more complicated than
in the case when all edges are rays emanating from one vertex.

In the case of several noncompact edges, it is natural to use the so-called
scattering data which are a generalization of the classical scattering data
for the Sturm-Liouville operator on the line. However, in the presence of
compact edges, the specification of the scattering data it is not sufficient
to determine the potential on all edges of the graph. Therefore, it is natural
to consider a mixture of the inverse spectral and inverse scattering problems.
Relevant data can be subdivided into two parts. The ''spectral'' part
consists of the so-called Weyl functions associated with the compact edges,
which generalize the Weyl function for the Sturm-Liouville operator on the
half-line and also on a noncompact tree with a single infinite edge [41].
Another part, in turn, includes a portion of negative eigenvalues together with
the so-called reflection coefficients and norming constants associated with all
but one noncompact edges. This part generalizes the left (or right) scattering
data for the Sturm-Liouville operator on the line. We note that in [32],
without formulations of results, a possible statement of the inverse problem
is discussed for star-type graph without compact edges.

Consider a noncompact star-type graph $\Gamma$ with the vertices
$v_0,\,v_1\ldots,v_p$ and the edges $\varepsilon_1,\ldots,\varepsilon_m,$ where
$v_0$ is an interior vertex; $\varepsilon_j=[v_0,v_j],\, j=\overline{1,p},$ are
compact edges, which are parameterized by the parameter $x_j\in[0,1],$ and
$\varepsilon_j,\,j=\overline{p+1,m}$ are infinite rays, which are parameterized by
the parameter $x_j\in[0,\infty).$ The value $x_j=0$ on all edges corresponds to
the point $v_0.$ Any function $y$ on $\Gamma$ can be represent in the form
$y=\{y_j\}_{j=\overline{1,m}},$ where the function $y_j(x_j)$ is defined
on $\varepsilon_j.$ For definiteness we suppose that $p\ge1$ and $m\ge p+2.$
Let $q=\{q_j\}_{j=\overline{1,m}}$ be a real-valued function on $\Gamma$;
it is called the potential. Assume that $q_j(x_j)\in L(0,1),\;j=\overline{1,p},$
and $q_j(x_j),\,x_jq_j(x_j)\in L(0,\infty),\;j=\overline{p+1,m}.$
Consider the Sturm-Liouville equation
$$
\ell_j y_j:=-y''_j+q_j(x_j)y_j=\la y_j, \quad j=\overline{1,m},             \eqno(5.1)
$$
where $x_j\in(0,1)$ for $j=\overline{1,p},$ and $x_j\in(0,\infty)$ for
$j=\overline{p+1,m}$. Let $y=\{y_j\}_{j=\overline{1,m}}$ satisfies the
standard matching conditions in the interior vertex $v_0:$
$$
y_1(0)=y_j(0),\; j=\overline{2,m},\qquad \sum_{j=1}^m y'_j(0)=0             \eqno(5.2)
$$
and the following boundary conditions in the boundary vertices:
$$
U_j(y_j):=y_j'(1)+H_jy(1)=0, \quad j=\overline{1,p}.                        \eqno(5.3)
$$
Boundary value problem (5.1)-(5.3) is denoted by $B.$ Let $\lambda=\rho^2.$
Denote ${\mathbb R}^*={\mathbb R}\setminus\{0\},$ $\Omega_\pm=\{\rho :\;
\pm Im\,\rho>0\}.$ Let $\psi_j(x_j,\la)$ be the solution of Eq. (5.1) on the
edge $\varepsilon_j$ for $j=\overline{1,p},$ satisfying initial conditions
$\psi_j(1,\lambda)=1,$ $U_j(\psi_j)=0.$ Let $e_j(x_j,\rho),\;\rho\in
\overline{\Omega_+}$ be the Jost solution (see [12]) of Eq. (5.1) on the
infinite edge $\varepsilon_j,$ $j=\overline{p+1,m}.$ Moreover, for
$j=\overline{1,m}$ we need the solutions $C_j(x_j,\lambda)$ and
$S_j(x_j,\lambda),$ satisfying the conditions
$C_j(0,\lambda)=S_j'(0,\lambda)=1,$ $C_j'(0,\lambda)=S_j(0,\lambda)=0.$

Let $\Psi_k(\lambda)=[\Psi_{kj}(x,\lambda)]_{j=\overline{1,m}},\,
k=\overline{1,p}$ be the solution of Eq. (5.1), satisfying matching
conditions (5.2) and the following boundary conditions:
$$
U_j(\Psi_{kj})=-\delta_{jk}, \quad j=\overline{1,p},
$$
$$
\Psi_{kj}(x_j,\lambda)=O(\exp(i\rho x_j)), \;\; x_j\to\infty, \;\;
\rho\in\overline{\Omega_+}, \;\; j=\overline{p+1,m}.
$$
The functions $\Psi_k$ and $M_k(\lambda):=\Psi_{kk}(1,\lambda)$ are called
{\it the Weyl solution} and {\it the Weyl function}, respectively, associated with
the boundary vertex $v_k,\,k=\overline{1,p}.$

Let the functions $f_k=\{f_{kj}\}_{j=\overline{1,m}},$ $k=\overline{p+1,m}$
be solutions of (5.1)-(5.3), and for $\rho\in\Omega_+$ satisfy the conditions
$$
f_{kk}(x_k,\rho)\sim\exp(-i\rho x_k),\; x_k\to\infty; \quad
f_{kj}(x_j,\rho)=O(\exp(i\rho x_j)),\;j=\overline{p+1,m}\setminus k,\;
x_j\to\infty.
$$
The function $f_k$ is called {\it the scattering solution}, associated with
the edge $\varepsilon_k$.

In order to construct Weyl solutions and scattering solutions one can use
auxiliary solutions $g_k=\{g_{kj}\}_{j=\overline{1,m}},$ $k=\overline{1,m},$
where
$$
g_{kj}(x_j,\rho)=e_j(x_j,\rho)\prod_{l\ne k,j}e_l(0,\rho), \quad j\ne k,
$$
$$
g_{kk}(x_k,\rho)=C_k(x_k,\lambda)\prod_{l\ne k}e_l(0,\rho)-
S_k(x_k,\lambda)\sum_{j\ne k}e_j'(0,\rho)\prod_{l\ne k,j} e_l(0,\rho),
$$
and $e_j(x_j,\rho)=\psi_j(x_j,\rho^2)$ for $j=\overline{1,p}.$

\smallskip
{\bf Theorem 5.1.} {\it For $j=\overline{1,m},$ the following representations hold
$$
\Psi_{kj}(x_j, \lambda)=\frac{1}{\Delta(\lambda)}g_{kj}(x_j, \rho),\;
k=\overline{1,p},\quad f_{kj}(x_j,\rho)=\frac{1}{a(\rho)}g_{kj}(x_j,\rho),
\;k=\overline{p+1,m}.
$$
where
$$
\Delta(\lambda)=\sum_{j=1}^m e_j'(0,\rho)\prod_{l\ne j} e_l(0,\rho), \quad
a(\rho)=\frac{1}{2i\rho}\Delta(\rho^2).
$$}

We note that Theorem 5.1 gives us the definition of scattering solutions
$f_k,$ $k=\overline{p+1,m}$ for $\rho\in{\mathbb R}^*_1:=
\{\rho:\rho\in{\mathbb R}^*,\, a(\rho)\ne0\}.$ One has
$$
f_{kk}(x_k,\rho)=e_k(x_k,-\rho)+s_k(\rho)e_k(x_k,\rho),
\quad \rho\in{\mathbb R}^*_1,\quad k=\overline{p+1,m}.
$$
The function $s_k(\rho)$ is called the {\it reflection coefficient},
associated with the edge $\varepsilon_j.$

The value of the parameter $\lambda,$ for which Eq. (5.1) has nonzero solutions
$y=[y_j(x)]_{j=\overline{1,m}},$ satisfying the matching conditions (5.2), the
boundary conditions (5.3) and the condition $y_j(x_j)\in L_2(0,\infty),$
$j=\overline{p+1,m},$ are called {\it eigenvalues} of $B,$ and the
corresponding solutions are called {\it eigenfunctions}.
If $y_j(x_j)\equiv 0$ for some $j$ with any eigenfunction $y,$ related to the
eigenvalue $\lambda_0,$ then $\lambda_0$ is called {\it invisible} from the edge
$\varepsilon_j;$ otherwise $\lambda_0$ is called {\it visible} from $\varepsilon_j.$

\smallskip
{\bf Theorem 5.2. } {\it Let $\Lambda$ be the set of visible eigenvalues from
noncompact edges $\varepsilon_j,$ $j=\overline{p+1,m-1}.$ Then this set is finite
and lies on the negative half-line:
$\Lambda=\{\lambda_n\}_{n=\overline{1,N}},\quad \lambda_n=\rho_n^2,\;
\rho_n=i\tau_n,\; \tau_n>0.$}

\smallskip
For brevity we put $f_k:=\Psi_k,$ $f_{kj}(x_j,\rho):=\Psi_{kj}(x_j,\rho^2),$
$j=\overline{1,m},$ $k=\overline{1,p}.$

\smallskip
{\bf Theorem 5.3.} {\it The function $f_k,\;k=\overline{1,m},$ has at most a simple
pole in the point $\rho_n,\,n=\overline{1,N}.$ The corresponding nonzero residue is
an eigenfunction related to the eigenvalue $\lambda_n=\rho_n^2.$}

\smallskip
Theorem 5.3 yields
$\mathop{\rm Res}\limits_{\rho=\rho_n}f_{kk}(x_k,\rho)=\alpha_{kn}e_k(x_k,\rho_n),$
where $\alpha_{kn}$ are called {\it norming numbers} for $\varepsilon_k.$ The set $J=\{M_\nu(\lambda),s_k(\rho),\lambda_n,\alpha_{kn}\}_{\nu=\overline{1,p},\,
k=\overline{p+1,m-1},\,n=\overline{1,N}}$ is called {\it the spectral data} for
the problem $B.$

\smallskip
{\bf Inverse problem 5.4.} Given $J,$ find $q=\{q_j\}_{j=\overline{1,m}}$
and $H_j,\;j=\overline{1,p}.$

\smallskip
The following theorem gives us the uniqueness of recovering $B$ from the spectral data.

\smallskip
{\bf Theorem 5.5. }{\it The specification of the spectral data $J$ uniquely
determines the potential $q$ and the coefficients $H_j$ $j=\overline{1,p}$.}

\smallskip
Theorems 5.1-5.5 are proved in [40]. A procedure for constructing the solution
jf inverse problem 5.4 is also obtained.

\medskip
{\bf 5.2. Inverse Sturm-Liouville problems on a noncompact sun-graph.}
We consider Sturm-Liouville differential operators on a sun-graph with one cycle and
an arbitrary number of noncompact edges-rays. We present results from [37] on the
inverse spectral problem for this important class of operators.

Let $G$ be a graph consisting of infinite rays $\mathcal{R}_k$,
$k=\overline{1,p},$ emanating from the vertices $v_k$, and a cycle $\mathcal{C}$,
going around these vertices. These vertices divide the cycle $\mathcal{C}$ into
the edges $\mathcal{E}_k$, $k=\overline{1,p}$. On the graph $G$ we consider
the differential expression
$$
\ell y=-y'' +q(x)y,
$$
where the derivative is taken with respect to the natural parameter (i.e. with
respect to the arch length on the edge containing the point $x$). This parameter
will be denoted by $|x|$. The potential $q(x)$ is real-valued and satisfies
the condition $(1+|x|)q(x)\in L(G).$ Denote by $\Lambda$ the set of eigenvalues
of the operator $L$, generated in $L_2(G)$  by the differential expression
$\ell y$ and the standard matching conditions in the vertices.

For each fixed $k=\overline{1,p},$ we define {\it the Weyl solution} $\psi_k(x,\rho),$
$x\in G,$ as a function on the graph with the properties:\\
1) $\ell \psi_k=\rho^2 \psi_k$ on each edge;\\
2) it satisfies the standard matching conditions in each vertex; \\
3) $\psi_k(x,\lambda)=O\left(\exp(i\rho|x|)\right)$ for $x\to\infty$,
$x\in\mathcal{R}_j$, $j=\overline{1,p}\setminus\{k\}$;\\
4) $\psi_k(x,\lambda)=\exp(-i\rho|x|)(1+o(1))$ for $x\to\infty$, $x\in\mathcal{R}_k$.

The function $\psi_k(x,\rho)$ is analytic in $\rho$ in the upper half-plane
$\Omega_+$ with the exception of a finite (possibly empty) set $Z_k^-$ of simple
poles situated on the imaginary axis. We note that if $\rho\in Z_k^-$, then
$\lambda=\rho^2\in\Lambda$, but the inverse assertion is not true. For the residues
${\rm res}_{\rho=\rho_0} \psi_k(x,\rho)$, $\rho_0\in Z_k^-$ we have
$$
{\rm res}_{\rho=\rho_0} \psi_k(x,\rho)=i\alpha_k(\rho_0)\exp(i\rho_0|x|)(1+o(1)),
\quad x\to\infty, x\in\mathcal{R}_k,
$$
where $\alpha_k(\rho_0)>0.$ Denote by $Z_0^+$ the set of
$\rho\in\mathbf{R}$ such that $\lambda=\rho^2\in\Lambda$. Then for all
$\rho_0\in\mathbf{R}\setminus\left(\{0\}\cup Z_0^+\right)$ there exist the limits
$\psi_k(x,\rho_0):=\lim\limits_{\rho\to\rho_0, \rho\in\Omega_+}\psi_k(x,\rho_0).$
If $\rho_0\in Z_0^+$, then $\psi_k(x,\rho)$ and $\psi'_k(x,\rho)$ are bounded for
$\rho\to\rho_0, \rho\in\Omega_+$. For the limit values $\psi_k(x,\rho)$,
$\rho\in\mathbf{R}\setminus\left(\{0\}\cup Z_0^+\right)$ the following relation holds
$$
\psi_k(x,\rho)=\exp(-i\rho |x|)+s_k(\rho)\exp(i\rho |x|)+o(1), \quad
x\to\infty, x\in\mathcal{R}_k.
$$
The function $s_k(\rho)$, $\rho\in\mathbf{R}\setminus (\{0\}\cup Z_0^+)$
is called {\it the reflection coefficient associated with $\mathcal{R}_k$}.
The set $J_k:=\{s_k(\rho),\rho\in\mathbf{R}\setminus (\{0\}\cup Z_0^+),
Z^-_k, \alpha_k(\rho), \rho\in Z_k^- \}$ is called {\it the scattering data
associated with $\mathcal{R}_k$}.

\smallskip
{\bf Theorem 5.6. }{\it The specification of the scattering data associated with
the ray $\mathcal{R}_k$, uniquely determines the potential on this ray.}

\smallskip
For recovering the potential on the cycle (and, consequently, on the whole graph $G$),
we need all sets $J_k$ associated with the rays, and additional information. For each
$k=\overline{1,p},$ we define {\it the additional Weyl solution}, which is defined
as a function $\psi_k^0(x,\rho)$, $x\in G$, $\rho\in\mathbf{R}\setminus\{0\}$ with
the properties:\\
1) $\psi_k^0(x,\rho)$ is a solution of equation $\ell y=\rho^2 y$;\\
2) $\psi_k^0(x,\rho)$ satisfies the standard matching conditions in the vertices;\\
3) the following asymptotic representation holds:
$$
\lim\limits_{x\to\infty,\, x\in\mathcal{R}_j}\exp(2i\rho|x|)\cdot
\frac{(\psi_k^0(x,\rho))'+i\rho\psi_k^0(x,\rho)}{(\psi_k^0(x,\rho))'
-i\rho\psi_k^0(x,\rho)}= -\frac{e'_j(v_j,-\rho)}{e'_j(v_j,\rho)},\;
j=\overline{1,p}\setminus \{k\},
$$
$$
\psi_k^0(x,\rho)=\exp(-i\rho|x|)+s_k^0(\rho)\exp(i\rho|x|)+o(1),
\; x\to\infty, x\in\mathcal{R}_k,
$$
where $e_j(x,\rho)$ is the classical Jost solution of the equation $\ell y=\rho^2 y$
on the ray $\mathcal{R}_j$. The function $s_k^0(\rho)$ is called {\it the additional
reflection coefficient associated with $\mathcal{R}_k$}.

Besides the scattering data and additional reflection coefficients associated
with the rays, we need spectral information related to directly the cycle
$\mathcal{C}$. Denote by $\Lambda^\mathcal{C}$ the spectrum of the problem
on the cycle (i.e. the spectrum of the periodic problem on $[0,\pi]$).
Let $\Lambda^k_\mathcal{C}$ be the spectrum of the problem on
the cycle with the Dirichlet boundary conditions at the vertex $v_k$, and
let $\Lambda_\mathcal{C}:=\bigcup\limits_{k=1}^p \Lambda_\mathcal{C}^k$.
One can show that $\Lambda^\mathcal{C}\setminus\Lambda_\mathcal{C}\subset\Lambda$.
On the set $\Lambda\cap (-\infty,0)=:\Lambda^-$ we define the function
$\omega(\lambda)$ as follows: $\omega(\lambda)=1,$ if
$\lambda\in\Lambda^\mathcal{C}\setminus\Lambda_\mathcal{C}$, and
$\omega(\lambda)=0$, otherwise. Let $\sigma^\mathcal{C}_n, n=\overline{1,\infty}$
be the sequence of signs for the elements from $\Lambda^p_\mathcal{C}$ associated
with the periodic problem on $\mathcal{C}$. The set
$$
J=\{J_k, k=\overline{1,p};\; s^0_k(\rho),\; \rho\in\mathbf{R}\setminus\{0\},
k=\overline{1,p}; \; \Lambda;\;  \omega(\lambda), \lambda\in\Lambda^-;
\; \sigma^\mathcal{C}_n, n=\overline{1,\infty}\}.
$$
is called {\it the complete scattering data}.

{\bf Theorem 5.7.} {\it The specification of the complete scattering data
uniquely determines the operator $L$.}

\smallskip
Theorems 5.6-5.7 are proved in [37]. A procedure for constructing the solution
of the inverse problem is also obtained.

\medskip
{\bf 5.3. Inverse Sturm-Liouville problems on a noncompact A-graph.}
We consider Sturm-Liouville differential operators on an arbitrary noncompact
A-graph. We remind that the A-graph has the property that any two cycles of
the graph can have at most one common point; in the rest the graph has an
arbitrary structure. This is a wide and important class of graphs. Inverse
problems for compact A-graphs were investigated in [17-18]. Noncompact case
produces essential qualitative modifications in the statement and in the
solution of inverse problems. In this subsection we briefly present the
results from [39] on inverse problems for noncompact A-graphs.

Let $G$ be a noncompact A-graph (see [18]) with the set of vertices
$V(G)$ and the set of edges $\mathcal{E}(G)\cup\mathcal{R}(G)$, where
$\mathcal{E}(G)$ is the set of compact edges, and $\mathcal{R}(G)$ is the
set of rays. Each edge is parameterized by the natural parameter $|x|$.
For two arbitrary points $x, x'$ on one edge, we denote by $|x-x'|$
the distance between these points along the edge (i.e. the length of the
corresponding arc). On $G$ we consider the differential expression
$$
\ell y:=-y''+q(x)y,                                                    \eqno(5.4)
$$
where $q$ is a real-valued function satisfying the condition
$(1+|x|)q(x)\in L(G).$ For each ray $r\in\mathcal{R}(G)$ we define
the scattering data $J_r$ associated with this ray analogously to the
case of the sun-graph (see subsection 5.2). For each boundary vertex $v$
we define the Weyl solution associated with $v$ as the function 
$\Phi_v(x,\lambda,G)$ with the properties:\\
1) $\ell \Phi_v=\lambda \Phi_v$, $x\in\mbox{int}\,r$, $r\in
\mathcal{E}(G)\cup\mathcal{R}(G)$;\\
2) it satisfies the standard matching conditions in all interior vertices
and the boundary conditions in all boundary vertices with the exception of
the vertex $v$;\\
3) $\Phi_v(\cdot,\lambda,G)\in L_2(G)$; $\Phi_v(v,\lambda,G)=1$.

The function $M_v(\lambda,G):=\partial_r\Phi_v(v,\lambda,G)$ is called
{\it the Weyl function associated with $v$}; here $\partial_r$ denotes
the derivative in the direction toward the interior of the edge $r$, which is
incident to $v.$

For each cycle $\mathfrak{c}$ we define the graph $G_\mathfrak{c}$ as follows.
Let $\mathfrak{c}$ consist of the edges $r_1, r_2,\ldots, r_p$, connecting
$v_0$ with $v_1$, $v_1$ with $v_2$, ..., $v_{p-1}$ with $v_0$, respectively, where
$v_0$ is the nearest vertex to the root, which belongs to $\mathfrak{c}$. Then
$G_\mathfrak{c}$ is the graph obtained from $G$ by the replacement of the edge
$r_p$ with the edge $r'_p$ of the same length connecting $v_{p-1}$ with an
additional vertex $v_\mathfrak{c}$. Put $q\left|_{r'_p}\right.=q\left|_{r_p}\right.$,
i.e. for arbitrary $x\in r'_p$ we define $q(x, G_\mathfrak{c}):=q(x',G)$,
where $x'\in r_p$ such that $|x'-v_{p-1}|=|x-v_{p-1}|$. On $G_\mathfrak{c}$
we consider the differential expression (5.4) together with the standard
matching conditions in the interior vertices and with the same boundary
conditions as for the graph $G$. The vertex $v_\mathfrak{c}$ is boundary
for $G_\mathfrak{c}$. For this vertex we introduce the Weyl function
$M_{v_\mathfrak{c}}(\lambda, G_\mathfrak{c})$.

The set $J=\{J_r, r\in\mathcal{R}, M_v(\cdot,G),v\in\partial G\setminus\{v^0\},
M_{v_\mathfrak{c}}(\cdot, G_\mathfrak{c}),\mathfrak{c}\in\mathcal{C}\}$ is
called the scattering data; here $\mathcal{C}$ is the set of cycles,
and $\partial G$ is the set of the boundary vertices of the graph $G.$

\smallskip
{\bf Theorem 5.8. }{\it The specification of the scattering data $J$
uniquely determines the potential $q$ on the graph $G.$}

\smallskip
Theorem 5.8 is proved in [39]. A procedure for constructing the solution
of the inverse problem is also obtained.

\medskip
{\bf 5.4. Variable order operators on a noncompact with a cycle.}
A number of practically important problems are reduced to variable order
differential equations on graphs, when orders of differential equations
on different edges can be different. Such problems are studied only recently,
and nowadays they are investigated insufficiently. In [48-50] inverse spectral
problems are considered for variable order operators on compact graphs
(see also subsections 4.2 and 4.3 of the present paper). In this subsection,
in contrast to the works [48-50], we study the inverse spectral problem on
a noncompact graph. The graph consists of a cycle and a ray with a common vertex.
On the ray we consider an equation of an arbitrary order $n>2.$
First we consider the case when the order of the equation on the cycle equals 2.

Let $\Gamma$ be a graph with the vertex $v_1$ and the edges $r_0$ and $r_1$,
where $r_1$ is the ray emanated from $v_1$, and $r_0$ is the cycle $[v_1,v_1]$
with the length $\pi.$ The edge $r_0$ is parameterized by the parameter
$x_0\in[0,\pi]$, and $r_1$ is parameterized by the parameter $x_1\in[0,\infty).$
A function $y$ on the graph can be represented as $y=(y_0(x_0), y_1(x_1)).$
On the edges $r_j,\; j=0,1,$ we consider the differential equations
$$
\ell_0 y_0:= -y_{0}''+q_{0}(x_0)y_0=\rho^2 y_0,                              \eqno(5.5)
$$
$$
\ell_1 y_1:= y_{1}^{(n)}+\sum_{j=0}^{n-2} q_{1j}(x_1)y_1^{(j)}=\rho^n y_1,   \eqno(5.6)
$$
where $q_{0}\in L(0,\pi)$ is a real-valued function, $n>2$, $q_{1j}$ are
complex-valued functions such that $q_{1j}(x_1)\exp(\tau x_1)\in L(0,\iy)$
for some $\tau>.0$ Let $\Omega_\nu=\{\rho: \arg\rho\in(\frac{\nu-1}{n}\pi,
\frac{\nu}{n}\pi)\},$ and let $\omega_k, k=\overline{1,n}$ be the roots
of the equation $\omega^n-1=0$ such that
$\mbox{Re}(\rho\omega_1)<\mbox{Re}(\rho\omega_2)<\ldots<\mbox{Re}(\rho\omega_n)$
for $\rho\in\Omega_\nu$. For $k=\overline{2,n},$ we define the Weyl solution
$\psi_k=\{\psi_{k0}(x_0,\rho), \psi_{k1}(x_1,\rho)\}$ as a function with
the properties:\\
1) $\psi_{k0}$ is a solution of Eq. (5.5), and $\psi_{k1}$ is a solution
of Eq. (5.6);\\
2) $\psi_k$ satisfies the matching conditions
$$
\psi_{k0}(0,\rho)=\psi_{k0}(\pi,\rho)=\psi_{k1}(0,\rho),\quad
\psi'_{k0}(\pi,\rho)=\psi'_{k0}(0,\rho)+\psi'_{k1}(0,\rho),
$$
$$
\psi^{(\nu-1)}_{k1}(0,\rho)=0, \nu=\overline{3,k};
$$
3) $\psi_{k1}(x_1,\rho)=\exp(\rho\omega_k x_1)(1+o(1)),\; x_1\to\infty$.

For $k=1$ we define $\psi_{11}(x_1,\rho)$ as the solution 0f Eq. (5.6)
such that\\
$\psi_{11}(x_1,\rho)=\exp(\rho\omega_1 x_1)(1+o(1)),$ $x_1\to\infty.$

\smallskip
{\bf Theorem 5.9.} {\it For $\nu=\overline{1,2n},\; k=\overline{1,n},$ the
function $\psi_{k1}(x_1,\rho)$ is holomorphic in $\Omega_\nu\setminus\Lambda_{k\nu}$
and continuous in $\overline{\Omega}_\nu\setminus\Lambda_{k\nu}$, where
$\Lambda_{k\nu}$ is an at most countable subset of $\overline{\Omega}_\nu$
without finite limit points. Moreover, $\psi_{k1}(x_1,\rho)=
O(\rho^{-N})$ for $\rho\to 0$, where $N$ is a positive number.}

We define the matrix $\Psi(x_1,\rho)$: $\Psi_{\nu k}:=\psi_{k1}^{(\nu-1)}$,
$k,\nu=\overline{1,n}$. Put $\Lambda_\nu:=\bigcup\limits_{k=1}^n\Lambda_{k\nu}$.

\smallskip
{\bf Theorem 5.10.} {\it There exists $\rho_*>0$ such that:\\
1) $\Lambda_{k\nu}\cap\Lambda_{m\nu}\cap\{|\rho|>\rho_*\}=\emptyset$;\\
2) for any $\rho_0\in\Lambda_\nu\cap\{|\rho|>\rho_*\}$ there exists a unique
strictly upper triangular matrix $\alpha(\rho_0)$ such that the function
$\Psi(x_1,\rho)(I-(\rho-\rho_0)^{-1}\alpha(\rho_0))$ is bounded in a neighborhood
of $\rho_0$.}

\smallskip
{\bf Theorem 5.11.} {\it For} $\rho\to\infty$, $\rho\in G_{\nu, \delta}:=
\{\rho\in\overline{\Omega}_\nu: \mbox{dist}(\rho,\Lambda_\nu)>\delta\}$, $\delta>0,$
$$
\Psi(x_1,\rho)=D_\rho\left[W\right]\exp(\rho x_1\omega)A(\rho), \ A(\rho)=O(1),
$$
{\it where} $D_\rho=\mbox{diag}(1, \rho,\ldots,\rho^{n-1})$, $\omega=\mbox{diag}
(\omega_1, \omega_2,\ldots,\omega_n)$, $W=(w_{jk})_{j,k=\overline{1,n}}$,
$w_{jk}=\omega_k^{j-1}$, $[W]=W+O(\rho^{-1})$, and $A(\rho)$ is {\it a upper triangular
matrix with 1 on the main diagonal.}

\smallskip
Let $\Sigma_\nu=\{\rho: \arg\rho=\frac{\nu}{n}\pi\}$. For $\rho_0\in\Sigma_\nu
\setminus\left(\Lambda_{\nu}\cup\Lambda_{\nu+1}\right)$ we put
$$
\Psi^-(x_1,\rho_0)=\lim\limits_{\rho\to\rho_0, \rho\in\Omega_\nu}\Psi(x_1,\rho),\
\Psi^+(x_1,\rho_0)=\lim\limits_{\rho\to\rho_0, \rho\in\Omega_{\nu+1}}\Psi(x_1,\rho).
$$
Then there exists a unique matrix $v(\rho_0)$ such that
$\Psi^+(x_1,\rho_0)=\Psi^-(x_1,\rho_0)v(\rho_0)$.

\smallskip
We will say that the collection of functions $\{q_0, q_{10},\ldots q_{1,n-2}\}$
belongs to the class $G$, if the assertion of Theorem 5.10 is valid for $\rho_*=0.$
In this case we define {\it the scattering data associated with the edge $r_1$}
as the collection
$J_1=\{v(\rho_0),\rho_0\in\Sigma_\nu\setminus (\Lambda_{\nu}\cup\Lambda_{\nu+1});
\Lambda_\nu; \alpha(\rho_0), \rho_0\in\Lambda_\nu,\; \nu=\overline{1,2n}\}.$

\smallskip
{\bf Theorem 5.12.} {\it The specification of the scattering data associated with
the edge $r_1$ uniquely determines the coefficients of Eq. (5.6).}

\smallskip
For recovering the coefficient of Eq. (5.5) we need additional data connected
with the periodic problem on $r_0$. Let $S_0(x_0,\la)$ and $C_0(x_0,\la)$
be solutions of the equation $\ell_0 y=\la y$ with the initial conditions $S_0(0,\la)=C'_0(0,\la)=0$, $S'_0(0,\la)=C_0(0,\la)=1.$ Put $F(\la):=C_0(\pi,\la)+S'_0(\pi,\la).$ Let $\{\lambda_m\}_{m=1}^\infty$ be the
spectrum of the boundary value problem for the equation $\ell_0 y=\lambda y$
with the Dirichlet boundary conditions at the endpoints of the interval
$(0,\pi),$ and let $\sigma_m\in\{-1,0,1\},\;m=\overline{1,\infty}$ be numbers
such that
$$
S'_{0}(\pi,\lambda_m)=F(\lambda_m)+\sigma_m\sqrt{F^2(\lambda_m)-4}/2.
$$
Denote by $N_0$ the set of all indices $m$ such that $\lambda_m$
is an eigenvalue of the periodic problem.

We define {\it the scattering data} as the set $J=\{J_1; \sigma_m,
m=\overline{1,\infty}; N_0; \lambda_m, m\in N_0\}$.

\smallskip
{\bf Theorem 5.13. } {\it The specification of the scattering data $J$
uniquely determines the coefficients of equations (5.5), (5.6).}

\smallskip
Theorems 5.9-5.13 are proved in [51]. A procedure for constructing the solution
of this inverse problem is also obtained.

\medskip
Now we consider essentially more difficult case when the order of the equation on
the cycle is more than two. The main difficulties here are that the periodic problem
on the cycle for higher order equations are essentially more difficult that for
the Sturm-Liouville operator. That is why the statement and the solution of the
inverse problem in this case are not trivial. Here we provide the results from
[52], where this complicated class of operators was studied in the first time.
For simplicity, we confine ourselves to the case when the order of the equation
on the cycle is equal to 3.

Let $\Gamma$ be a graph consisting of the cycle $r_0$ of length $T$ and the
ray $r_1$ emanated from some point $v_1\in r_0$. A function $y$ on $\Gamma$
may be represented in the form $y=\{y_0(x), x\in[0,T],\;y_1(x),x\in[0,\iy)\}$.
On the cycle $r_0$ we consider the equation
$$
\ell_0 y_0\equiv
D^3 y_0+p_{01}(x) Dy_0+p_{00}(x)y_0=\rho^3 y_0,\quad D=-id/dx,              \eqno(5.7)
$$
where $\ell_0^*=\ell_0$. On the ray $r_1$ we consider the equation
$$
\ell_1 y_1\equiv D^N y_1+\sum_{s=0}^{N-2} p_{1s}(x) D^s y_1=\rho^N y_1,    \eqno(5.8)
$$
where $N\geq 3$ and $|p_{1s}(x)|\exp(\tau x)\in L(0,\iy)$ for some
$\tau>0.$ Consider the linear forms
$$
U_\nu(y):=\sigma_\nu y^{(\nu-1)}(0)+\sum_{s=0}^{\nu-2}\sigma_{\nu s}y^{(s)}(0),
u_{\xi\nu}(y)=(-1)^{\chi_{\xi\nu}}y^{(\nu-1)}(\xi),\; \nu=1,2,
$$
where $\xi\in\{0,T\}$, $\chi_{0\nu}=0$, $\chi_{T\nu}=\chi$, $\chi_{T3}=\chi+1$,
$\chi\in\{0,1\}$. For a function $y=(y_0, y_1)$ and $\nu\in\overline{1,N}$ we define
the matching condition $C(\nu)$ by the relation $u_{0\nu}(y_0)=u_{T\nu}(y_0)=
U_\nu(y_1)$, and the condition $K(\nu)$ by the relation $u_{0\nu}(y_0)+
u_{T\nu}(y_0)+U_\nu(y_1)=0$ for $\nu\leq 3,$ and $U_\nu(y_1)=0$ for $\nu>3.$

Let $S_l:=\{\arg (i\rho)\in((l-1)\frac{\pi}{N},l\frac{\pi}{N}\}$. For
a fixed $l$ we denote by $R_k$, $k=\overline{1,N}$ the roots of the equation
$R^N-1=0$ such that $\mbox{Re}(i\rho R_1)<
\mbox{Re}(i\rho R_2)<\ldots<\mbox{Re}(i\rho R_N)$ for all $\rho\in S_l$.
Fix $\chi\in\{0,1\}$. For each $k=\overline{1,N}$, in each sector $S_l$
we define the Weyl solution $\psi_k=\{\psi_{k0}(x,\rho),
\psi_{k1}(x,\rho)\}$ of order $k$ as the solution of system (5.7)-(5.8)
with the condition $\psi_{k1}(x,\rho)=\exp(i\rho R_k x)(1+o(1)),\,
x\to\infty$ and with the matching conditions $C(\nu)$, $\nu=\overline{1,\nu_k-1},
\nu_k=\min\{k,3\}$, $K(\nu)$, $\nu=\overline{\nu_k,k}$. The solution
$\psi_k$ exists and is unique for all $\rho\in\bar S_l$ with the exception
of an at most countable set without finite limit points (besides, possibly,
the point 0).

{\it Condition $G_0$.} For each $k$ the function $\psi_{k1}(x,\rho)$ is
holomorphic in $S_l\cap\{|\rho|<\delta\}$ for some $\delta>0$, continuous in
$\bar S_l\cap\{|\rho|<\delta\}\setminus\{0\}$, and
$\psi^{(\nu-1)}_{k1}(x,\rho)=O\left(\rho^{-M}\right),\; k,\nu=\overline{1,N}$
for $\rho\to 0$, where $M<\infty$.

Let $Y^\alpha_{k1}(x,\rho)$, $k=\overline{1,N}$ be the fundamental system of
solutions (in each sector $S_l$) of Eq. (5.8) constructed in [11] and having
the following properties:\\
1) the functions $Y^\alpha_{k1}(x,\rho)$ are holomorphic in $\rho$ in
$S_l\cap\{|\rho|>\rho_\alpha\}$, and $\rho_\alpha\to 0$ for $\alpha\to\infty$;\\
2) $\lim\limits_{x\to\infty} Y^\alpha_{k1}(x,\rho)\exp(-i\rho R_k x)=1$;\\
3) $D^\nu Y^{\alpha}_{k1}(x,\rho)=(\rho R_k)^\nu\exp(i\rho R_k x)[1]$,
$[1]:=\left(1+O(\rho^{-1})\right)$, $x\geq\alpha$, $\rho\to\infty$.

Fix $\alpha$ such that the functions $Y^\alpha_{k1}(x,\rho)$ are holomorphic in
$S_l\cap\{|\rho|>\delta/2\}$ and continuous in $\bar S_l\cap\{|\rho|\geq\delta/2\}$.
Then for $\rho\in \bar S_l\cap\{|\rho|\geq\delta/2\}$, the following
representations hold:
$$
\psi_{k1}(x,\rho)=Y^\alpha_{k1}(x,\rho)
+\sum_{j<k}\gamma^{\alpha}_{jk}(\rho)Y^\alpha_{j1}(x,\rho),\;
\psi_{k0}(x,\rho)=\sum_{j=1}^3 \beta_{jk}(\rho) C_j(x,\lambda),             \eqno(5.9)
$$
where $\lambda=\rho^3$, and $C_j(x,\lambda)$ are solutions of the equation
$\ell_0 y=\lambda y$ under the conditions $C_j^{(\nu-1)}(0,\lambda)=\delta_{j\nu}$.

Substituting (5.9) into the matching conditions, we obtain a system of linear
algebraic equa\-tions with respect to $\gamma^{\al}_{jk}(\rho)$, $\beta_{jk}(\rho)$.
Denote by $\Delta_k(\rho)$ the determinant of this system, and by $Z_{kl}$
the set of its zeros situated in $\bar S_l\setminus\{0\}$ for $k\in\overline{2,N}$,
and in $\bar S_l$ for $k=1.$ Clearly, the functions $\psi_k$ are continuous on
$\bar S_l\setminus\left(\{0\}\cup Z_{kl}\right)$, and are holomorphic in
$S_l\setminus Z_{kl}$. The functions $Y^\alpha_{k1}(x,\rho)$ admit an
analytic continuation to a domain of the form
$S^\varepsilon_l\setminus\{|\rho|\leq\rho_\alpha\}$, where
$S^\varepsilon_l:=S_l+i\varepsilon\exp(i(l-\frac{1}{2})\frac{\pi}{N})$.
Therefore, for any $\rho_0\in Z_{kl}$ the functions $\psi_{kj}(x,\rho)$
admit a holomorphic continuation in a deleted neighborhood of $\rho_0$.
In the sequel, we will assume that the following condition is fulfilled.

{\it Condition $G_1$.} The sets $Z_{kl}$ for different $k$ do not intersect.
Each $\rho_0\in Z_{kl}$ is a simple zero of $\Delta_k(\rho)$ and a simple
zero of $\psi_{kj}(x,\rho)$ (at least for one $j\in\{0,1\}$).

For any $\rho_0\in Z_{kl}$ there exist numbers $v^l_{jk}(\rho_0)$,
$j<k$ such that the functions
$$
D^{\nu-1}\psi_{k1}(x,\rho)-(\rho-\rho_0)^{-1}
\sum_{j<k}v^l_{jk}(\rho_0)D^{\nu-1}\psi_{j1}(x,\rho_0),\; \nu=\overline{1,N}
$$
are holomorphic in a neighbourhood of $\rho_0$.

For studying the behavior of the Weyl solutions for large $\rho$,
we rewrite them in the form:
$$
\psi_{k1}(x,\rho)=Y^0_{k1}(x,\rho)+\sum_{j<k}\gamma_{jk}(\rho)
Y^0_{j1}(x,\rho),\quad
\psi_{k0}(x,\rho)=\sum_{j=1}^3\beta_{jk}(\rho)Y_{j0}(x,\rho),          \eqno(5.10)
$$
where $Y_{j0}(x,\rho)$ are Birchhoff's solutions [55]. Substituting (5.10)
into the matching conditions, we obtain a system of linear algebraic
equations with respect to $\gamma_{jk}(\rho)$, $\beta_{jk}(\rho)$.
Solving this system we find
$$
\gamma_{sk}(\rho)=-D_{sk}(\rho)/D_k(\rho),
$$
where $D_k(\rho)$ is the determinant of the system, and $D_{sk}(\rho)$
can be obtained from $D_k(\rho)$ by formal replacement $Y_{s1}$ with
$Y_{k1}$. Now we consider the following system (it is obtained from
the above-mentioned system by the replacement $Y_{kj}$ with their main
asymptotical terms):
$$
\sum_{s=1}^3\beta_{ks}\omega^{s(\nu-1)}=(-1)^{\chi_{T\nu}}\sum_{s=1}^3\beta_{ks}
\omega^{s(\nu-1)}\exp(i\rho\omega^s T)=\sigma_\nu\left(\sum_{s=1}^{k-1}
\gamma_{ks}R_s^{\nu-1}+R_k^{\nu-1}\right),\; \nu=\overline{1,\nu_k-1},
$$
$$
\sum_{s=1}^3\beta_{ks}\omega^{s(\nu-1)}+(-1)^{\chi_{T\nu}}\sum_{s=1}^3\beta_{ks}
\omega^{s(\nu-1)}\exp(i\rho\omega^s T)+\sigma_\nu\left(\sum_{s=1}^{k-1}
\gamma_{ks}R_s^{\nu-1}+R_k^{\nu-1}\right)=0,\;  \nu=\overline{\nu_k,3},
$$
$$
\sum_{s=1}^{k-1}\gamma_{ks}R_s^{\nu-1}+R_k^{\nu-1}=0,\; \nu=\overline{4,k}.
$$
Denote by $D^0_k(\rho)$ the determinant of this system. One has
$$
D^0_k(\rho)=A_{k0}+\sum_{m=1}^3 (A^+_{km}
\exp(i\rho\omega^m T)+A^-_{km}\exp(-i\rho\omega^m T),
$$
where the numbers $A_{k0}$, $A^\pm_{km}$ depend only on the coefficients
$\sigma_\nu$ of the forms $U_\nu$ and the sector $S_l$. For the determinates
$D_k(\rho)$, $D_{sk}(\rho)$ the following asymptotical formulae hold
$$
D_k(\rho)=[A_{k0}]+\sum_{m=1}^3 ([A^+_{km}]
\exp(i\rho\omega^m T)+[A^-_{km}]\exp(-i\rho\omega^m T)),
$$
$$
D_{sk}(\rho)=[A_{sk0}]+\sum_{m=1}^3 ([A^+_{skm}]
\exp(i\rho\omega^m T)+[A^-_{skm}]\exp(-i\rho\omega^m T)),
$$
where the numbers $A_{sk0}$, $A^\pm_{skm}$ also depend only on the coefficients
$\sigma_\nu$ of the forms $U_\nu$ and the sector $S_l$. Let the regularity
condition $A^\pm_{km}\neq 0$ is fulfilled for all $k,m,l$.

Fix $\chi\in\{0,1\}$. We define the matrix
$\Psi=(\Psi_{\nu k})_{k,\nu=\overline{1,N}}$ via
$\Psi_{\nu k}(x,\rho):=D^{\nu-1}\psi_{k1}(x,\rho)$. For
$\rho_0\in\Sigma_l\setminus(Z_l\cup Z_{l+1})$ (where $\Sigma_l:=
\bar S_l\cap\bar S_{l+1}$, $Z_l:=\bigcup_k Z_{kl}$) we define
$\Psi_-(x_1,\rho_0):=\lim\limits_{\rho\to\rho_0, \rho\in S_{l}}
\Psi(x,\rho)$, $\Psi_+(x,\rho_0):=\lim_{\rho\to\rho_0, \rho\in
S_{l+1}}\Psi(x,\rho)$ and the matrix $v(\rho_0):=\Psi^{-1}_-(x,\rho_0)
\Psi_+(x,\rho_0)$. Furthermore, for $\rho_0\in Z_{kl}$ we define the matrices
$v_l(\rho_0):=\left(v^l_{jk}(\rho_0)\right)_{j,k=\overline{1,N}}$.

The set
$$
J_1^\chi=\left\{v(\rho), \rho\in \Sigma_l\setminus (Z_l\cup Z_{l+1}), Z_{kl},
v_l(\rho), \rho\in Z_{kl}, k=\overline{1,N}, l=\overline{1,2N}\right\}
$$
is called the scattering data associated with the ray $r_1$.

{\bf Theorem 5.14. }{\it The specification of the scattering data $J^\chi_1$
uniquely determines the potentials $p_{1s}$, $s=\overline{0,N-2}$ and $\psi$.}

\smallskip
Denote by $\Lambda^\pm$ the spectra of the boundary value problems for the equation
$\ell_0 y=\lambda y$ with the conditions $y^{(\nu-1)}(0)\pm y^{(\nu-1)}(T)=0$, $\nu=\overline{1,3}$. Denote by $\Lambda^\pm_s$, $s=1,2$ the spectra of the
problems for the same equation with the conditions
$y^{(s-1)}(0)=y^{(s-1)}(T)=y^{(2-s)}(0)\pm y^{(2-s)}(T)=0.$ Let
$\Delta^\pm_{30}(\lambda)$ and $\Delta^\pm_{3s}(\lambda)$ be the characteristic
functions of these problems, respectively. Put $\Lambda^\pm_{3s}:=
\Lambda^\pm\cap\Lambda^\pm_s$.

The set $J=\{J_1^0, J_1^1,\Lambda^+_{31},\Lambda^-_{31},\Lambda^+_{32},\Lambda^-_{32}\}$
is called {\it the global scattering data}.

\smallskip
{\bf Theorem 5.15.} {\it The specification of the global scattering data $J$
uniquely determines the potentials $p_{1s}$, $s=\overline{0,N-2}$ and
$p_{0s}$, $s=0,1$.}

\smallskip
Theorems 5.14-5.15 are proved in [52]. A constructive procedure for the solution
of this inverse problem is also obtained.

\medskip
{\bf Acknowledgment.} This work was supported by Grant 1.1436.2014K of the Russian
Ministry of Education and Science and by Grant 13-01-00134 of Russian Foundation for
Basic Research.

\begin{center}
{\bf Bibliography}
\end{center}

\begin{enumerate}
\item[{[1]}] Montrol E. Quantum theory on a network. J. Math. Phys.
     11, no.2 (1970), 635-648.
\item[{[2]}] Nicaise S. Some results on spectral theory over networks,
     applied to nerve impulse transmis\-sion. Lect. Notes Math., vol. 1771,
     Springer, 1985, 532-541.
\item[{[3]}] Langese J.E., Leugering G. and Schmidt J.P.G. Modelling,
     analysis and control of dynamic elastic multi-link structures.
     Birkh\"auser, Boston, 1994.
\item[{[4]}] Kottos T. and Smilansky U. Quantum chaos on graphs.
     Phys. Rev. Lett. 79 (1997), 4794-4797.
\item[{[5]}] Tautz J., Lindauer M. and Sandeman D.C. Transmission
     of vibration across honeycombs and its detection by bee leg receptors.
     J. Exp. Biol. 199 (1999), 2585-2594.
\item[{[6]}] Dekoninck B. and Nicaise S. The eigenvalue problem for
     networks of beams. Linear Algebra Appl. 314 (2000), no. 1-3, 165-189.
\item[{[7]}] Sobolev A. and Solomyak M., Schr\"odinger operator on
     homogeneous metric trees: spectrum in gaps. Rev. Math. Phys. 14,
     no.5 (2002), 421-467.
*\item[{[8]}] Pokorny Yu.V., Penkin O.M., Pryadiev V.L., Borovskikh A.V.,
     Lazarev K.P., Shabrov C.A. Differential equations on geometrical graphs.
     Moscow, Fizmatlit, 2004.
\item[{[9]}] Marchenko V.A. Sturm-Liouville operators and their applications.
     "Naukova Dumka",  Kiev, 1977;  English  transl., Birkh\"auser, 1986.
\item[{[10]}] Levitan B.M. Inverse Sturm-Liouville problems. Nauka, Moscow,
     1984; English transl., VNU Sci.Press, Utrecht, 1987.
\item[{[11]}] Yurko V.A. Method of Spectral Mappings in the Inverse
     Problem Theory, Inverse and Ill-posed Problems Series. VSP, Utrecht, 2002.
\item[{[12]}] Yurko V.A. Introduction to the theory of inverse spectral
     problems. Moscow, Fizmatlit, 2007, 384pp. (Russian).
\item[{[13]}] Belishev M.I. Boundary spectral Inverse Problem on a class of graphs
     (trees) by the BC method. Inverse Problems 20 (2004), 647-672.
\item[{[14]}] Yurko V.A. Inverse spectral problems for Sturm-Liouville operators
     on graphs. Inverse Problems, 21, no.3 (2005), 1075-1086.
\item[{[15]}] Brown B.M.; Weikard R. A Borg-Levinson theorem for trees. Proc. R.
     Soc. Lond. Ser. A Math. Phys. Eng. Sci. 461 (2005), no.2062, 3231-3243.
\item[{[16]}] Yurko V.A. Inverse problems for Sturm-Liouville operators on
     bush-type graphs. Inverse Problems, 25, no.10 (2009), 105008, 14pp.
\item[{[17]}] Yurko V.A. An inverse problem for Sturm-Liouville operators
     on A-graphs.  Applied Math. Letters, 23, no.8 (2010), 875-879.
\item[{[18]}] Yurko V.A. Reconstruction of Sturm-Liouville differential operators
     on A-graphs. Differ. Uravneniya, 47, no.1 (2011), 50-59 (Russian); English
     transl. in Differential Equations, 47, no.1 (2011), 50-59.
\item[{[19]}] Yurko V.A. An inverse problem for Sturm-Liouville operators on
     arbitrary compact spatial networks. Doklady Akad. Nauk. 432, no.3 (2010),
     318-321 (Russian); English transl. in Doklady Mathematics 81, no.3 (2010),
     410-413.
\item[{[20]}] Yurko V.A. Inverse spectral problems for differential operators
     on arbitrary compact graphs. Journal of Inverse and Ill-Posed Proplems,
     18, no.3 (2010), 245-261.
\item[{[21]}] Freiling G. and Yurko V.A. Inverse problems for differential
     operators on graphs with general matching conditions. Applicable Analysis,
     86, no.6 (2007), 653-667.
\item[{[22]}] Avdonin S.; Kurasov P. Inverse problems for quantum trees.
     Inverse Problems and Imaging, vol.2, no.1 (2008), 1-21.
\item[{[23]}] Carlson R. Inverse eigenvalue problems on directed graphs.
     Trans. Amer. Math. Soc. 351, no.10 (1999), 4069-4088.
\item[{[24]}] Freiling G., Ignatiev M., Yurko V.A. An inverse spectral problem
     for Sturm-Liouville operators with singular potentials on star-type graphs.
     Analysis on graphs and its applica\-tions. Proceedings of Symposia in Pure
     Mathematics, vol.77, 397-408, Amer. Math. Soc., Providence, 2008.
\item[{[25]}] Kurasov P. Inverse scattering for lasso graph. J. Math. Phys.
     54 (2013), no.4, 042103, 14pp.
\item[{[26]}] Pivovarchik V.N. Inverse problem for the Sturm-Liouville equation
     on a simple graph. SIAM J. Math. Anal. 32, no.4 (2000), 801-819.
\item[{[27]}] Yang C-Fu; Yang X-P. Uniqueness theorems from partial information
     of the potential on a graph. Journal of Inverse and Ill-Posed Problems, 19,
     no. 4-5 (2011), 631-639.
\item[{[28]}] Yurko V.A. Inverse nodal problems for differential operators on graphs.
     Journal of Inverse and Ill-Posed Problems, 16, no.7 (2008), 715-722.
\item[{[29]}] Freiling G. and Yurko V.A. Inverse nodal problems for differential
     operators on graphs with a cycle. Tamkang Journal of Mathematics, 41, no.1
     (2010), 15-24.
\item[{[30]}] Yurko V.A. Recovering differential pencils on compact graphs.
     Journal of Diff. Equations, 244, no.2 (2008), 431-443.
\item[{[31]}] Yurko V.A. An inverse problem for differential pencils on graphs with
    a cycle. Journal of Inverse and Ill-Posed Problems 22, no.5 (2014), 625-641.
\item[{[32]}] Gerasimenko N.I. Inverse scattering problem on a noncompact graph.
     Teoret. Mat. Fiz. 74 (1988), no. 2, 187-200 (Russian); English transl. in
     Theor. Math. Phys. 75 (1988), 460-470.
\item[{[33]}] Harmer M. Inverse scattering for the matrix Schr\"odinger operator
     and Schr\"odinger operator on graphs with general self-adjoint boundary
     conditions. ANZIAM J. 44 (2002), 161-168.
\item[{[34]}] Kurasov P.; Stenberg F. On the inverse scattering problem on
     branching graphs. J. Phys. A: Math. Gen. 35 (2002), 101-121.
\item[{[35]}] Pivovarchik V.; Latushkin Y. Scattering in a forked-shaped waveguide.
     Integral Equat. Oerator Theory 61 (2008), no.3, 365-399.
\item[{[36]}] Trooshin I.; Marchenko V.; Mochizuki K. Inverse scattering on
     a graph containing circle. Analytic methods of analysis and differ.
     equations: AMADE 2006, 237-243, Camb. Sci. Publ., Cambridge, 2008.
\item[{[37]}] Ignatyev M.Yu. and Freiling G. Spectral analysis for the Sturm-Liouville
     operator on sun-type graphs. Inverse Problems 27, no.9 (2011), 095003 (17pp).
\item[{[38]}] Ignatyev M.Yu. Inverse scattering problem for Sturm-Liouville operator
     on one-vertex noncompact graph with a cycle. Tamkang J. Math. 42, no.3 (2011),
     365-384.
\item[{[39]}] Ignatyev M.Yu. Inverse scattering problem for Sturm-Liouville operator
     on non-compact A-graph. Uniqueness result. arXiv:1311.2862 (2014)
\item[{[40]}] Buterin S.A. and Freiling G. Inverse spectral-scattering problem for
     the Sturm-Liouville operator on a noncompact star-type graph, Tamkang Journal
     of Mathematics 44 (2013), no.3, 327-349.
\item[{[41]}] Freiling G. and Yurko V.A. Inverse spectral problems for
     Sturm-Liouville operators on noncompact trees. Results in Mathematics,
     50, no.3-4 (2007), 195-212.
\item[{[42]}] Yurko V.A. An inverse spectral problem for pencils of differential
     operators on noncompact networks. Differ. Uravneniya, 44, no.12 (2008),
     1658-1666 (Russian); English transl. in Differential Equations, 44, no.12
     (2008), 1721-1729.
\item[{[43]}] Yurko V.A. Inverse problems for Bessel-type differential equations
     on noncompact graphs using spectral data. Inverse Problems, 27, no.4
     (2011), 045002, 17pp.
\item[{[44]}] Yurko V.A. Inverse problems for differential operators of any order
     on trees. Matem. Zametki, 83, no.1 (2008), 139-152 (Russian); English transl.
     in Mathematical Notes, 83, no.1 (2008), 125-137.
\item[{[45]}] Yurko V.A. Inverse spectral problems for arbitrary order differential
     operators on noncompact trees. Journal of Inverse and Ill-Posed Problems,
     20, no.1 (2012), 111-132.
\item[{[46]}] Yurko V.A. Inverse problems for differential systems on graphs
     with regular singularities. Math. Notes  96, no.4 (2014), 617-621.
\item[{[47]}] Pokornyi Yu.V., Beloglazova T.V., Lazarev K. On a class of variable
     order ordinary differen\-tial equation on a graph. Matem. Zametki 73, no.3
     (2003), 469-470; English transl. in Mathem. Notes 73, no.3 (2003).
\item[{[48]}] Yurko V.A. Recovering variable order differential operators on
     star-type graphs from spectra. Differ. Uravneniya, 49, no.12 (2013),
     1537-1548. (Russian); English transl. in  Differ. Equa\-tions 49, no.12
     (2013), 1490-1501.
\item[{[49]}] Yurko V.A. Inverse problems for differential operators of variable
     orders on star-type graphs: general case. Analysis and Mathematical Physics,
     vol. 4, no.3 (2014), 247-262.
\item[{[50]}] Bondarenko N.P. Inverse problems  for the differential operator on
     the graph with a cycle with different orders on different edges. Tamkang
     Journal of Mathematics (accepted).
\item[{[51]}] Ignatyev M.Yu. Uniqueness of recovering a variable order differential
     operators on the simplest noncompact graph with a cycle. Mathematics and
     Mechanics, vol.15, Saratov University Press, Saratov, 2013, 35-37 (Russian).
\item[{[52]}] Ignatyev M.Yu. Uniqueness of the solution of the inverse scattering
     problem for a variable order differential operators on the simplest
     noncompact graph with a cycle. Izvest. Saratov University. Ser. Math.,
     Mech., Inform. vol. 14, no.4, part 2 (2014), 542-549 (Russian).
\item[{[53]}] Currie S. and Watson B. The M-matrix for the Sturm-Liouville
     equation on graphs. Proc. Royal Soc. Edinburgh 139A (2009), 775-796.
\item[{[54]}] Yurko V.A. An inverse problems for the matrix Sturm-Liouville
     equation on a finite interval. Inverse Problems, 22, no.4 (2006), 1139-1149.
\item[{[55]}] Naimark M.A. Linear differential operators. 2nd ed.,
     Nauka, Moscow, 1969; English transl. of 1st ed., Parts I,II, Ungar,
     New York, 1967, 1968.
\end{enumerate}


\begin{tabular}{ll}
Name:             &   Yurko, Vjacheslav  \\
Place of work:    &   Department of Mathematics, Saratov State University \\
{}                &   Astrakhanskaya 83, Saratov 410012, Russia \\
E-mail:           &   yurkova@info.sgu.ru\\
\end{tabular}

\end{document}